\documentclass[11pt]{amsart}
\usepackage{amssymb}
\usepackage[all]{xy}

\usepackage{mathrsfs}

\newtheorem{thm}{Theorem}[section]
\newtheorem{prop}[thm]{Proposition}
\newtheorem{defin}[thm]{Definition}
\newtheorem{lemma}[thm]{Lemma}
\newtheorem{cor}[thm]{Corollary}

\theoremstyle{remark}
\newtheorem{remark}[thm]{Remark}

\numberwithin{equation}{section}

\newcommand{\nc}{\newcommand} \nc{\on}{\operatorname}

\nc{\pa}{\partial}

\nc{\g}{{\mathfrak g }}
\nc{\cA}{{\mathcal A}} \nc{\cB}{{\cal B}}\nc{\cC}{{\mathcal C}} 
\nc{\cD}{{\mathcal D}} 
\nc{\cE}{{\mathcal E}} \nc{\cG}{{\mathcal G}}\nc{\cH}{{\cal H}} 
\nc{\cI}{{\cal I}} \nc{\cJ}{{\cal J}}\nc{\cK}{{\cal K}} 
\nc{\cL}{{\mathcal L}} \nc{\cR}{{\cal R}} \nc{\cS}{{\mathcal S}}   
\nc{\cV}{{\cal V}} \nc{\cX}{{\cal X}}

\nc{\Maj}{\on{Maj}}
\nc{\sh}{\on{sh}}
\nc{\Id}{\on{Id}}\nc{\Diff}{\on{Diff}}
\nc{\Perm}{\on{Perm}}\nc{\conc}{\on{conc}}\nc{\Alt}{\on{Alt}}
\nc{\ad}{\on{ad}}\nc{\Der}{\on{Der}}\nc{\End}{\on{End}}
\nc{\no}{\on{no\ }} \nc{\res}{\on{res}}\nc{\ddiv}{\on{div}}
\nc{\Sh}{\on{Sh}} \nc{\card}{\on{card}}\nc{\dimm}{\on{dim}}
\nc{\Sym}{\on{Sym}} \nc{\Jac}{\on{Jac}}\nc{\Ker}{\on{Ker}}
\nc{\Spec}{\on{Spec}}\nc{\Cl}{\on{Cl}}
\nc{\Imm}{\on{Im}}\nc{\limm}{\lim}\nc{\Ad}{\on{Ad}}
\nc{\ev}{\on{ev}} \nc{\Hol}{\on{Hol}}\nc{\Det}{\on{Det}}
\nc{\Bun}{\on{Bun}}\nc{\diag}{\on{diag}}\nc{\pr}{\on{pr}} 
\nc{\Span}{\on{Span}}\nc{\Comp}{\on{Comp}}\nc{\Part}{\on{Part}}
\nc{\tensor}{\on{tensor}}\nc{\ind}{\on{ind}}\nc{\id}{\on{id}}
\nc{\Hom}{\on{Hom}}\nc{\Quant}{\on{Quant}} \nc{\Dequant}{\on{Dequant}}\nc{\Def}{\on{Def}}
\nc{\AutLBA}{\on{AutLBA}}\nc{\AutQUE}{\on{AutQUE}}
\nc{\LBA}{{\on{LBA}}}\nc{\Aut}{\on{Aut}}\nc{\QUE}{{\on{QUE}}}
\nc{\Lyn}{\on{Lyn}}\nc{\Cof}{\on{Cof}}\nc{\LCA}{\underline{\on{LCA}}}\nc{\FLBA}{\on{FLBA}}
\nc{\LA}{\on{LA}}\nc{\FLA}{\on{FLA}}\nc{\EK}{\on{EK}}
\nc{\class}{\on{class}}\nc{\br}{\on{br}}\nc{\co}{\on{co}}
\nc{\Prim}{\on{Prim}}\nc{\ren}{\on{ren}}\nc{\lbr}{\on{lbr}}
\nc{\SC}{\on{SC}}\nc{\ol}{\overline}\nc{\FL}{\on{FL}}
\nc{\FA}{\on{FA}}\nc{\alg}{{\on{alg}}}\nc{\KZ}{\on{KZ}}
\nc{\op}{{\on{op}}}\nc{\cop}{{\on{cop}}}
\nc{\Inn}{\on{Inn}}\nc{\OutDer}{\on{OutDer}}
\nc{\inv}{{\on{inv}}}\nc{\gr}{{\on{gr}}}\nc{\Lie}{{\on{Lie}}}
\nc{\Out}{{\on{Out}}}\nc{\univ}{{\on{univ}}}
\nc{\GT}{{\on{GT}}}\nc{\GRT}{{\on{GRT}}}
\nc{\GS}{{\on{GS}}} 
\nc{\restr}{{\on{restr}}}\nc{\lin}{{\on{lin}}}
\nc{\mult}{{\on{mult}}}\nc{\qt}{{\on{qt}}}
\nc{\compl}{{\on{compl}}} \nc{\Trees}{{\on{Trees}}}
\nc{\Ord}{{\on{Ord}}} \nc{\norm}{{\on{norm\ ord}}}
\nc{\fd}{{\on{fd}}}\nc{\Maps}{{\on{Maps}}}
\nc{\CYBE}{{\on{CYBE}}} \nc{\CY}{{\on{C}}} \nc{\can}{{\on{can}}}
\nc{\Alg}{{\underline{\on{Alg}}}} \nc{\Poisson}{{\underline{\on{Poisson}}}}
\nc{\Coalg}{{\underline{\on{Coalg}}}} \nc{\Vect}{{\underline{\on{Vect}}}}
\nc{\UE}{{\on{UE}}}\nc{\Bialg}{{\underline{\on{Bialg}}}}
\nc{\QTA}{{\underline{\on{QTA}}}}\nc{\QTQUE}{{\underline{\on{QTQUE}}}}
\nc{\QTLBA}{{\underline{\on{QTLBA}}}}\nc{\TA}{{\underline{\on{TA}}}}
\nc{\MT}{{\underline{\on{MT}}}}\nc{\QYBE}{{\underline{\on{QYBE}}}}
\nc{\QTBialg}{{\underline{\on{QTBialg}}}}
\nc{\qcocomm}{{\on{quasi-cocomm}}}
\nc{\qcomm}{{\on{quasi-comm}}}\nc{\cocomm}{{\on{cocomm}}}
\nc{\dual}{{\on{dual}}} \nc{\sym}{{\on{sym}}}
\nc{\cycl}{{\on{cycl}}} \nc{\Prop}{{\underline{\on{Prop}}}}
\nc{\uV}{{\underline{V}}} \nc{\uR}{{\underline{R}}}
\nc{\Mor}{{\on{Mor}}}\nc{\Ob}{{\on{Ob}}}
\nc{\Bil}{{\on{Bil}}} \nc{\invt}{{\on{invt}}}
\nc{\udim}{{\operatorname{\underline{dim}}}}
\nc{\ddim}{{\operatorname{\underline{\underline{dim}}}}}

\nc{\al}{\alpha}\nc{\de}{\delta}
\nc{\eps}{\epsilon}\nc{\la}{{\lambda}}
\nc{\si}{\sigma}\nc{\z}{\zeta}

\nc{\La}{\Lambda}

\nc{\ve}{\varepsilon} \nc{\vp}{\varphi} 

\nc{\AAA}{{\mathbb A}}\nc{\BB}{{\mathbb B}}
\nc{\CC}{{\mathbb C}}\nc{\ZZ}{{\mathbb Z}} 
\nc{\QQ}{{\mathbb Q}} \nc{\NN}{{\mathbb N}}\nc{\VV}{{\mathbb V}} 
\nc{\KK}{{\mathbb K}} 
\nc{\LL}{{\mathfrak{Lib}}} 

\nc{\ff}{{\mathbf f}}\nc{\bg}{{\mathbf g}}
\nc{\ii}{{\mathbf i}}\nc{\kk}{{\mathbf k}}
\nc{\bl}{{\mathbf l}}\nc{\zz}{{\mathbf z}} 
\nc{\pp}{{\mathbf p}}\nc{\qq}{{\mathbf q}} 

\nc{\cF}{{\cal F}}\nc{\cM}{{\cal M}}\nc{\cO}{{\cal O}}
\nc{\cT}{{\cal T}}\nc{\cW}{{\cal W}}\nc{\cP}{{\mathcal P}}

\nc{\Assoc}{{\mathbf Assoc}}
\nc{\ul}{\underline}

\nc{\ub}{{\underline{b}}}
\nc{\uk}{{\underline{k}}} 
\nc{\un}{{\underline{n}}} \nc{\um}{{\underline{m}}}
\nc{\up}{{\underline{p}}}\nc{\uq}{{\underline{q}}}
\nc{\ur}{{\underline{r}}}
\nc{\us}{{\underline{s}}}\nc{\ut}{{\underline{t}}}
\nc{\uw}{{\underline{w}}}
\nc{\uz}{{\underline{z}}}
\nc{\ual}{{\underline{\alpha}}}\nc{\ualpha}{{\underline{\alpha}}}
\nc{\ubeta}{{\underline{\beta}}}\nc{\ugamma}{{\underline{\gamma}}}
\nc{\ueps}{{\underline{\epsilon}}}\nc{\ueta}{{\underline{\eta}}}
\nc{\uzeta}{{\underline{\zeta}}}\nc{\ula}{{\underline{\lambda}}}
\nc{\umu}{{\underline{\mu}}}\nc{\unu}{{\underline{\nu}}}
\nc{\usigma}{{\underline{\sigma}}}\nc{\utau}{{\underline{\tau}}}
\nc{\uI}{{\underline{I}}}\nc{\uJ}{{\underline{J}}}
\nc{\uK}{{\underline{K}}}\nc{\uM}{{\underline{M}}}
\nc{\uN}{{\underline{N}}}

\nc{\A}{{\mathfrak a}}
\nc{\G}{{\mathfrak g}}
\nc{\x}{{\mathfrak x}}
\nc{\B}{{\mathfrak b}} \nc{\C}{{\mathfrak c}} 
\nc{\D}{{\mathfrak d}} \nc{\HH}{{\mathfrak h}}
\nc{\iii}{{\mathfrak i}}\nc{\mm}{{\mathfrak m}}\nc{\N}{{\mathfrak n}} 
\nc{\ttt}{{\mathfrak{t}}}\nc{\U}{{\mathfrak u}}\nc{\V}{{\mathfrak v}}
\nc{\grt}{{\mathfrak{grt}}}\nc{\gt}{{\mathfrak gt}}
\nc{\SL}{{\mathfrak{sl}}}\nc{\out}{{\mathfrak{out}}}

\nc{\SG}{{\mathfrak S}}

\nc{\wt}{\widetilde} \nc{\wh}{\widehat}

\nc{\bn}{\begin{equation}}\nc{\en}{\end{equation}} \nc{\td}{\tilde}

\def\C{{\mathbb{C}}}

\def\Lib{{\mathfrak{Lib}}}
\def\<<{\langle\langle}
\def\>>{\rangle\rangle}

\def\mt{{\mathfrak{mt}}}
\def\dmr{{\mathfrak{dmr}}}

\author{Benjamin Enriquez}
\author{Hidekazu Furusho}

\address{Institut de Recherche Math\'{e}matique Avanc\'{e}e, UMR 7501, 
Universit\'{e} de Strasbourg et CNRS, 7
rue Ren\'{e} Descartes, 67000 Strasbourg, France}
\email{enriquez@math.unistra.fr}

\address{Graduate School of Mathematics, Nagoya University, 
Furo-cho, Chikusa-ku, Nagoya, 464-8602, Japan}
\email{furusho@math.nagoya-u.ac.jp}

\date{November 24, 2017}
\title{A stabilizer interpretation of double shuffle Lie algebras}

\begin{document}

\bibliographystyle{amsalpha+}
\maketitle

\begin{abstract}
According to Racinet's work,
the scheme of double shuffle and regularization relations between cyclotomic analogues 
of multiple zeta values has the structure of a torsor over a pro-unipotent $\mathbb Q$-algebraic group $\sf{DMR}_0$, which is an 
algebraic subgroup of a pro-unipotent $\mathbb Q$-algebraic group of outer automorphisms of a free Lie algebra. We show that the 
harmonic (stuffle) coproduct of double shuffle theory may be viewed as an element of a module over the above group, and that 
$\sf{DMR}_0$ identifies with the stabilizer of this element. 
We identify the tangent space at origin of $\sf{DMR}_0$ with the stabilizer Lie algebra of the harmonic coproduct, thereby obtaining 
an alternative proof of Racinet's result stating that this space is a Lie algebra (the double shuffle Lie algebra).
\end{abstract}

%
%
%
%
%
%
%
%
%
%
%
%
%

\tableofcontents

\section{Introduction}
Multiple $L$-values (MLV in short) $L(k_1,\cdots,k_m;\zeta_1,\cdots,\zeta_m)$
are the complex numbers defined by the following series
\begin{equation*}\label{multiple L-value}
L(k_1,\cdots,k_m;\zeta_1,\cdots,\zeta_m)
:=\sum_{n_1>\cdots>n_m>0}\frac{\zeta_1^{n_1}\cdots \zeta_m^{n_m}}
{n_1^{k_1}\cdots n_m^{k_m}}
\end{equation*}
for $m$, $k_1$,\dots, $k_m\in{\mathbb Z}_{>0}$
and $\zeta_1$,\dots,$\zeta_m$ in the group $\mu_N$ of $N$-th roots of unity in $\mathbb C$, where $N$ is an integer $\geq 1$.
They converge if and only if $(k_1,\zeta_1)\neq (1,1)$. Multiple zeta values are regarded as a special case for
$N=1$. These values have recently garnered much interest due to their appearance in various fields of physics and mathematics
(\cite{BK,Dr,LM}). 
In connection with motive theory (\cite{De,DG}),
linear and algebraic relations among MLV's are particularly important.
The extended (regularized) double shuffle relations (\cite{IKZ, Rac}) 
might be one of the most fascinating ones.

Racinet's formalism \cite{Rac} is quite useful to state these relations. The basic ingredients of this formalism are: 
a noncommutative formal power series algebra 
$\mathbb C\langle\langle X\rangle\rangle$ over variables $x_0,x_\sigma$ 
($\sigma\in\mu_N$), another noncommutative formal power series algebra $\mathbb C\langle\langle Y\rangle\rangle$ over another set of variables, a map 
$\mathbb C\langle\langle X\rangle\rangle^\times\to\mathbb C\langle\langle Y\rangle\rangle^\times$ denoted $G\mapsto G_\star$
(see Def.\ \ref{def:3:6}), 
a coproduct map
$\Delta:\mathbb C\langle\langle X\rangle\rangle\to\mathbb C\langle\langle X\rangle\rangle^{\hat\otimes 2}$
for which  the elements of $X$ are primitive
and another coproduct map
$\Delta_*:\mathbb C\langle\langle Y\rangle\rangle\to\mathbb C\langle\langle Y\rangle\rangle^{\hat\otimes 2}$ (see \S\ref{sect:hccpo}). 
The central object of Racinet's approach is a certain invertible non-commutative formal power series
$\mathscr{I}$ in $\C\<< X\>>^\times$, constructed through iterated integrals and whose coefficients are expressed in terms of MLV's;
in particular, the coefficient of $x_0^{k_1-1}x_{\zeta_1}x_0^{k_2-1}x_{\zeta_1\zeta_2}\cdots x_0^{k_m-1}x_{\zeta_1\cdots\zeta_m}$
is 
$$L(k_1,\cdots,k_m;\zeta_1,\dots,\zeta_m).
$$
The regularized double shuffle relations for MLV's can then be formulated as
\begin{equation}\label{eqs:sat:by:I}
\Delta(\mathscr{I})=\mathscr{I}\otimes \mathscr{I} \text{ and }
\Delta_*(\mathscr{I}_\star)=\mathscr{I}_\star\otimes \mathscr{I}_\star
\end{equation}
(equalities in $\C\<< X\>>^{\hat\otimes 2}$, resp., $\C\<< Y\>>^{\hat\otimes 2}$). 
The series $\mathscr{I}$ has constant term equal to 1, coefficients of $x_0$ and $x_1$ equal to 0, and coefficient of $x_0x_1$ 
equal to $-(2\pi\sqrt{-1})^2/24$. 
For $\Bbbk$ a commutative $\mathbb Q$-algebra and $\lambda\in\Bbbk$, the set of series in 
$\Bbbk\langle\langle X\rangle\rangle$ satisfying (\ref{eqs:sat:by:I}) and the above coefficient conditions, with $2\pi\sqrt{-1}$ 
replaced by $\lambda$, is denoted
by ${\sf{DMR}}_\lambda(\Bbbk)$. So $\mathscr{I}$ belongs to ${\sf{DMR}}_{2\pi\sqrt{-1}}(\mathbb C)$. 
The structure of the regularized double shuffle relations is clarified by the following result: 

\begin{thm} \label{thm:racinet} (\cite{Rac}  Th.\ I, \S 3.2.3)
The scheme $\Bbbk\mapsto{\sf{DMR}}_0(\Bbbk)$ forms a pro-unipotent subgroup scheme of $\sf{MT}$. For 
$\Bbbk$ a $\mathbb Q$-commutative algebra and $\lambda\in\Bbbk$, the set ${\sf{DMR}}_\lambda(\Bbbk)$
is a torsor (principal homogeneous space) over the group ${\sf{DMR}}_0(\Bbbk)$.
\end{thm}

Here $\sf{MT}$ is\footnote{${\sf{DMR}}$ and $\sf{MT}$ stand for the French ``double m\'elange et r\'egularisation'' and ``groupe de Magnus 
tordu''.} the pro-unipotent $\mathbb Q$-algebraic group scheme 
$$
\Bbbk\mapsto{\sf{MT}}(\Bbbk):=(\{\text{series in }\Bbbk\langle\langle X\rangle\rangle\text{ with 
constant term equal to }1\},\circledast),
$$ 
where the group structure $\circledast$ is such that there is a group morphism from ${\sf{MT}}(\Bbbk)$ to the group of 
automorphisms of the topological $\Bbbk$-algebra $\Bbbk\langle\langle X\rangle\rangle$ taking $x_0$ to itself and each 
$x_\sigma$, $\sigma\in\mu_N$, to a conjugate of itself, and which are invariant with respect to the natural action of 
$\mu_N$ (see \S\ref{sect:qgsatla}). 

Thm.\ \ref{thm:racinet} shows that the scheme ${\sf{DMR}}_{2\pi\sqrt{-1}}(\mathbb C)$, which corresponds to the double shuffle 
relations of MLVs, is isomorphic to the group ${\sf{DMR}}_0(\mathbb C)$. This group is prounipotent, and therefore isomorphic
to its Lie algebra, which is the degree completion of a positively graded Lie algebra $\mathfrak{dmr}_0$. It is therefore an important 
question, which is still open today, to determine the size of this
Lie algebra (more precisely, its Poincar\'e polynomial, as $\mathfrak{dmr}_0$ is graded). 

In order to better understand $\mathfrak{dmr}_0$, it seems useful to improve our understanding of its Lie algebra structure. 
Here it seems helpful to recall the situation with the Grothendieck-Teichm\"uller Lie algebra $\mathfrak{grt}_1$. This Lie 
algebra was introduced by Drinfeld (\cite{Dr}) and later shown to be contained, up to the automorphism $(x_0,x_1)\mapsto(x_0,-x_1)$, 
in $\mathfrak{dmr}_0$ when $N=1$ (\cite{Fur}). A transparent 
interpretation of its Lie algebra structure was given by Ihara (\cite{I}): $\mathfrak{grt}_1$ is isomorphic to any of the Lie algebras
of symmetric special outer derivations of $\mathfrak{p}_n$ for $n\geq 5$, where $\mathfrak p_n$ is the graded Lie algebra
associated to the prounipotent completion of the pure sphere braid group with $n$ strands. 

It is also instructive to recall the situation with the Kashiwara-Vergne Lie algebra $\mathfrak{krv}_2$ arising in the work of 
Alekseev and Torossian (\cite{AT}). The inclusion $\mathfrak{grt}_1\subset\mathfrak{krv}_2$ was proved in 
\cite{AT}, while the inclusion $\mathfrak{dmr}_0\subset \mathfrak{krv}_2$ was proved in \cite{Sch}. In contrast to the previous 
situation, the Lie algebra status of $\mathfrak{krv}_2$ is rather transparent from its definition, as it is related to a Lie algebra cocycle. 

Let us review the known proofs of the Lie algebra nature of $\mathfrak{dmr}_0$. The space $\mathfrak{dmr}_0$ is introduced in 
\cite{Rac} as the tangent space of ${\sf{DMR}}_0$, viewed as a subspace of $\sf MT$, at the unit element. The proof of
\cite{Rac} that $\mathfrak{dmr}_0$ is a Lie algebra appeals to coderivations of the coalgebra 
$(\mathbb C\langle\langle Y\rangle\rangle,\Delta_*)$;  in the case when $N=1$,
this proof was streamlined in the appendix of \cite{Fur}. The statement that $\dmr_0$ is a Lie algebra was also announced without proof in 
\cite{Ec}; in the case when $N=1$, it was later proved using the techniques of Ecalle's theory of moulds in \cite{SaSch}.

In this paper, we introduce a graded Lie subalgebra $(\Lib(X),\langle,\rangle)$ of a graded Lie algebra $\mt$, whose completion  
is the Lie algebra of ${\sf MT}$, and show that {\it $\mathfrak{dmr}_0$ coincides, up to addition of an abelian Lie algebra, with the stabilizer 
of an element in a $(\Lib(X),\langle,\rangle)$-module (see Thm. \ref{main:thm});} this gives an alternative proof of the Lie algebra nature of 
$\mathfrak{dmr}_0$.  
We then prove a group version of this result: for any commutative $\mathbb Q$-algebra $\Bbbk$, we construct a subgroup 
$(\mathrm{exp}(\hat{\Lib}_{\Bbbk}(X)),\circledast)$ of ${\sf MT}(\Bbbk)$ and prove that {\it ${\sf DMR}_0(\Bbbk)$ coincides (up to a 
central additive group) with the stabilizer of an element of an $(\mathrm{exp}(\hat{\Lib}_{\Bbbk}(X)),\circledast)$-module.} 
More precisely: 

 \par{1)} there is a group representation of ${\sf MT}(\Bbbk)$ on the topological 
$\Bbbk$-module $\Bbbk\langle\langle Y\rangle\rangle$; we denote by $g\mapsto S^Y_g$ the corresponding group morphism 
${\sf MT}(\Bbbk)\to\mathrm{Aut}_{\Bbbk}^{\mathrm{cont}}(\Bbbk\langle\langle Y\rangle\rangle)$ (see (\ref{module:kY}));  

\par{2)} there is a group morphism 
$\Theta_\Bbbk:(\mathrm{exp}(\hat{\Lib}_{\Bbbk}(X)),\circledast)\to{\sf MT}(\Bbbk)$ (see \S\ref{sect:cot}); this morphism
should be viewed as a correction of the canonical inclusion related with regularization;  

\par{3)} the set $\widetilde{{\sf DMR}}_0(\Bbbk):=\{e^{\beta x_1}\cdot g\cdot e^{\alpha x_0}\ |\ \alpha,\beta\in\Bbbk, \ 
g\in{\sf DMR}_0(\Bbbk)\}$ (where $\cdot$ is the usual product in $\Bbbk\langle\langle X\rangle\rangle^\times$) is a subgroup of 
$(\mathrm{exp}(\hat{\Lib}_{\Bbbk}(X)),\circledast)$, isomorphic to ${\sf DMR}_0(\Bbbk)\times\Bbbk^2$ (see \S\ref{tgsaimd});   

\noindent 
we then prove (see \S\ref{sect:iofktd}): 
\begin{thm}\label{thm:tsdcwtsoted}
The subgroup $\widetilde{{\sf DMR}}_0(\Bbbk)$ coincides with the stabilizer of the element $\Delta_*$ of the space
$\mathrm{Hom}^{\mathrm{cont}}_{\Bbbk}(\Bbbk\langle\langle Y\rangle\rangle,\Bbbk\langle\langle Y\rangle\rangle^{\hat\otimes 2})$ 
of all continuous $\Bbbk$-linear maps $\Bbbk\langle\langle Y\rangle\rangle\to\Bbbk\langle\langle Y\rangle\rangle^{\hat\otimes 2}$, 
equipped with the action of $(\mathrm{exp}(\hat{\LL}(X)),\circledast)$, pulled back by $\Theta_\Bbbk$ of the natural action of 
${\sf MT}(\Bbbk)$, namely
\begin{align*}\widetilde{{\sf DMR}}_0(\Bbbk)=
 \{g\in\mathrm{exp}(\hat{\Lib}_\Bbbk(X))
\ |\ (S^Y_{\Theta_\Bbbk(g)})^{\otimes 2}\circ\Delta_*=\Delta_*\circ S^Y_{\Theta_\Bbbk(g)}
(\text{in }\mathrm{Hom}^{\mathrm{cont}}_{\Bbbk}(\Bbbk\langle\langle Y\rangle\rangle,\Bbbk\langle\langle Y\rangle\rangle^{\hat\otimes 2})
)\}. 
\end{align*}
\end{thm} 
Therefore $\widetilde{{\sf DMR}}_0(\Bbbk)$ may be interpreted as the set of $g\in\mathrm{exp}(\hat{\Lib}_\Bbbk(X))$ such that the  
diagram $\xymatrix{ \Bbbk\langle\langle Y\rangle\rangle\ar^{\Delta_*}[r]\ar_{S^Y_{\Theta_\Bbbk(g)}}[d] & \Bbbk\langle\langle Y\rangle\rangle^{\hat\otimes 2}\ar^{(S^Y_{\Theta_\Bbbk(g)})^{\otimes 2}}[d]\\ 
\Bbbk\langle\langle Y\rangle\rangle \ar^{\Delta_*}[r]& \Bbbk\langle\langle Y\rangle\rangle^{\hat\otimes 2} }$ is commutative, i.e. such that $S^Y_{\Theta_\Bbbk(g)}$ is an automorphism of the 
topological coalgebra $(\Bbbk\langle\langle Y\rangle\rangle,\Delta_*)$. 

The realization of this program is done in the following steps. \S\ref{section:1} presents background material necessary for the proof of the 
Lie algebra version of Thm. \ref{thm:tsdcwtsoted}, which is obtained in \S\ref{sect:cola} (Thm. \ref{main:thm}). In \S\ref{sect:eotlaht}, 
we construct the group homomorphism  
$\Theta_\Bbbk$, a necessary step for the proof of Thm. \ref{thm:tsdcwtsoted} in \S\ref{sect:coqgs}.

\section{Lie algebras of derivations of free Lie algebras}\label{section:1}

In this section, we present the material needed for the comparison of $\mathfrak{dmr}_0$ and $\mathfrak{stab}(\Delta_*)$.  
In \S \ref{sect:1:1} and \S\ref{sect:hccpo}, we recall the formalism of double shuffle theory: the Ihara bracket on 
the augmentation ideal $\mathbb Q\langle X\rangle_0$ of the associative algebra $\mathbb Q\langle X\rangle$ and 
on the free Lie algebra $\LL(X)$, the maps $\pi_Y,\mathrm{corr}:
\mathbb Q\langle X\rangle\to\mathbb Q\langle Y\rangle$, the endomorphism $\mathbf q$ of $\mathbb Q\langle Y\rangle$, 
and the harmonic coproduct $\Delta_*$ of $\mathbb Q\langle Y\rangle$. In \S \ref{sect:1:2:new}, we construct modules 
$\mathbb Q\langle X\rangle,\mathbb Q\langle Y\rangle$ over the space $\mathbb Q\langle X\rangle_0$, viewed as a Lie algebra 
via its Ihara bracket, and a morphism $\mathbb Q\langle X\rangle\stackrel{\mathbf q\circ\pi_Y}{\to}\mathbb Q\langle Y\rangle$ 
between these modules. In \S \ref{sect:fctalah}  we construct a Lie algebra morphism $\theta:\LL(X)\to\mathbb Q\langle X\rangle_0$, 
where both sides are equipped with the Lie algebra structures arising from the Ihara bracket (Prop. \ref{tmtltq}). One then pulls back the module 
$\mathbb Q\langle Y\rangle$ over $\mathbb Q\langle X\rangle_0$ to get a module over $\LL(X)$, and constructs from there a new 
module $\mathrm{Hom}_{\mathbb Q}(\mathbb Q\langle Y\rangle,\mathbb Q\langle Y\rangle^{\otimes 2})$ over this Lie algebra 
containing a particular element $\Delta_*$. In \S \ref{sect:stabilizer}, we construct the stabilizer subalgebra $\mathfrak{stab}(\Delta_*)\subset\LL(X)$ of 
this element. Finally, in \S \ref{sect:relations}, we show a relation between $\theta$, a morphism 
$\mathrm{sec}:\mathbb Q\langle Y\rangle\to\mathbb Q\langle X\rangle$ and an endomorphism $\tilde{\mathbf p}$ of 
$\mathbb Q\langle X\rangle$. 

\subsection{Ihara brackets and outer derivation Lie algebras}\label{sect:1:1}

\subsubsection{Ihara brackets}

Let $\Gamma$ be a finite commutative group, whose product will be denoted multiplicatively. Let $X$ be the alphabet 
$\{x_0\}\sqcup \{x_\gamma|\gamma\in\Gamma\}$, indexed by 
$\{0\}\sqcup\Gamma$. Let $\mathbb Q\langle X\rangle$ be the free noncommutative associative algebra with unit over the alphabet $X$ and 
$\LL(X)$ be the Lie subalgebra of  $\mathbb Q\langle X\rangle$ generated by $X$ (see \cite{Rac}, \S 2.2.3).

The group $\Gamma$ acts on $X$ by $\sigma\cdot x_0=x_0$, $\sigma\cdot x_\gamma=x_{\sigma\gamma}$
for $\gamma\in\Gamma$. This action extends to an action on $\mathbb Q\langle X\rangle$ and $\LL(X)$, which will be 
denoted $\sigma\mapsto t_\sigma$ (see \cite{Rac}, \S 3.1.1). 

For $\psi\in\mathbb Q\langle X\rangle$, 
let $d_\psi$ be the derivation of $\mathbb Q\langle X\rangle$ given by 
$$
\forall \sigma\in\Gamma,\quad d_\psi(x_0)=0,\quad d_\psi(x_\sigma)=[x_\sigma,t_\sigma(\psi)]
$$
(see \cite{Rac}, \S 3.1.12.2) and let $s_\psi$ be the linear endomorphism of $\mathbb Q\langle X\rangle$ defined 
by 
$$
s_\psi(v):=\psi v+d_\psi(v)
$$ for $v\in\mathbb Q\langle X\rangle$ (see \cite{Rac}, \S 3.1.12.1). 

Let $\mathbb Q\langle X\rangle_0$ be the subspace of $\mathbb Q\langle X\rangle$ of all series with zero constant term
(this is denoted $\mathfrak{mt}(k)$ in \cite{Rac}, \S 3.1.8). On $\mathbb Q\langle X\rangle_0$, one defines the (Ihara)
bracket 
\begin{equation}\label{Ihara bracket}
\langle\psi_1,\psi_2\rangle:=s_{\psi_1}(\psi_2)-s_{\psi_2}(\psi_1) 
\end{equation}
(see \cite{Rac}, (3.1.10.2)). This restricts to a bracket on the subspace $\LL(X)$ of $\mathbb Q\langle X\rangle_0$. We have therefore 
a Lie algebra inclusion 
\begin{equation}\label{LA:incl}
(\LL(X),\langle,\rangle)\subset(\mathbb Q\langle X\rangle_0,\langle,\rangle)=:\mt. 
\end{equation}

\subsubsection{Outer derivation Lie algebras}\label{sect:odla}

Call an associative (resp., Lie) algebra derivation of $\mathbb Q\langle X\rangle$ (resp., $\LL(X)$) 
{\it tangential} if it takes any of the generators $x_\gamma$, $\gamma\in\Gamma\cup\{0\}$
to an element of $[x_\gamma,\mathbb Q\langle X\rangle]$ (resp., $[x_\gamma,\LL(X)]$). The tangential 
derivations of either kind form a Lie algebra, of which the subspace of inner derivations 
is an ideal. In each case, the quotient algebra has a natural action of $\Gamma$. The subalgebra of invariants of this action 
is denoted $\mathfrak{outder}^*(\mathbb Q\langle X\rangle)^\Gamma$ (resp., $\mathfrak{outder}^*(\LL(X))^\Gamma$). 
As $\mathbb Q\langle X\rangle$ is the universal enveloping algebra 
of $\Lib(X)$, any derivation of the latter algebra extends uniquely to a derivation of the former algebra, which induces an injective 
Lie algebra morphism $\mathfrak{outder}^*(\Lib(X))^\Gamma\hookrightarrow \mathfrak{outder}^*(\mathbb Q\langle X\rangle)^\Gamma$. 

Introduce a $\mathbb Z^2$-grading in the algebras $\LL(X)$ and $\mathbb Q\langle X\rangle$ by $\mathrm{deg}(x_0):=(1,0)$, 
$\mathrm{deg}(x_\gamma):=(0,1)$ for $\gamma\in\Gamma$. This grading is compatible with the Lie bracket $\langle,\rangle$ in 
both cases. Moreover, it induces a $\mathbb Z^2$-grading on $\mathfrak{outder}^*(\LL(X))^\Gamma$ and 
$\mathfrak{outder}^*(\mathbb Q\langle X\rangle)^\Gamma$.    

\begin{lemma}\label{lemma:taiozgla}
There are isomorphisms of $\mathbb Z^2$-graded Lie algebras
\begin{equation}\label{eq:Lib}
(\LL(X),\langle,\rangle)\simeq \mathfrak{outder}^*(\LL(X))^\Gamma\oplus (\mathbb Q x_0\oplus\mathbb Q x_1), 
\end{equation}
and 
\begin{equation}\label{eq:mt}
\mt\simeq \mathfrak{outder}^*(\mathbb Q\langle X\rangle)^\Gamma\oplus (\mathbb Q[x_0]_0\oplus\mathbb Q[x_1]_0), 
\end{equation}
in which $\mathbb Q[x_0]_0,\mathbb Q[x_1]_0$ are the spaces of polynomials in $x_0,x_1$ with no constant terms, and 
where the spaces $\mathbb Q x_0\oplus\mathbb Q x_1$ and $\mathbb Q[x_0]_0\oplus\mathbb Q[x_1]_0$  are equipped with 
abelian Lie algebra structures. The Lie algebra inclusion (\ref{LA:incl}) is compatible with the isomorphisms (\ref{eq:Lib}) and (\ref{eq:mt}) 
and with the natural inclusions of the respective summands of the right-hand side of (\ref{eq:Lib}) in the summands of the right-hand side 
of (\ref{eq:mt}). 
\end{lemma}

\proof

The map $\psi\mapsto(\text{class of }d_\psi)$ induces a surjective Lie algebra morphism 
$(\mathbb Q\langle X\rangle_0,\langle,\rangle)
\to\mathfrak{outder}^*(\mathbb Q\langle X\rangle)^\Gamma$. 
Its kernel is $\mathbb Q[x_0]_0\oplus\mathbb Q[x_1]_0$, which is contained in the center of 
$(\mathbb Q\langle X\rangle_0,\langle,\rangle)$.

For $n\geq 0$, the components of $\mathbb Q[x_0]_0$ and of $\mathbb Q\langle X\rangle_0$ of 
bidegree $(n,0)$ are the same, so $\mathfrak{outder}^*(\mathbb Q\langle X\rangle)^\Gamma 
[n,0]=0$ for $n\geq 0$ (where $X[n,m]$ is the subspace of a bigraded space $X$ with bidegree 
$(n,m)$). 

The direct sum of homogenous components of $\mathbb Q\langle X\rangle$ of bidegrees $(0,n)$, $n\geq 0$ identifies with 
the free associative algebra $\mathbb Q\langle x_\gamma,\gamma\in\Gamma\rangle$. Let $I$ be the two-sided ideal of this 
algebra generated by the family $\{x_\gamma\ |\ \gamma\in\Gamma\setminus\{1\}\}$. This is a $\mathbb Q$-linear subspace of 
$\mathbb Q\langle X\rangle$, complementary to $\mathbb Q[x_1]$, so $I\oplus\mathbb Q[x_1]_0=\oplus_{n\geq 0}
\mathbb Q\langle X\rangle_0[0,n]$. It follows that $\mathfrak g :=I\oplus(\oplus_{n,m>0}\mathbb Q\langle X\rangle_0[n,m])$ 
is such that $\mathfrak g\oplus(\mathbb Q[x_0]_0\oplus\mathbb Q[x_1]_0)=\mathbb Q\langle X\rangle_0$.
One also has $\langle I,I\rangle\subset I$, so that $\mathfrak g$ is a Lie subalgebra of $\mathbb Q\langle X\rangle_0$. 
This proves both that the latter equality is a direct sum decomposition of the Lie algebra $(\mathbb Q\langle X\rangle_0,\langle,\rangle)$, 
and that the natural projection 
$\mathfrak g\to\mathfrak{outder}^*(\mathbb Q\langle X\rangle)^\Gamma$ is an isomorphism, which proves (\ref{eq:mt}).    

The proof of (\ref{eq:Lib}) is similar, the inclusion $I\subset\mathbb Q\langle X\rangle_0$ 
being replaced by the inclusion in the Lie algebra $\LL(x_\gamma,\gamma\in\Gamma)$ of its 
ideal generated by $\{x_\gamma\ |\ \gamma\in\Gamma\setminus\{1\}\}$. \hfill \qed\medskip 

\subsection{Harmonic coproduct, correction and projection operators}\label{sect:hccpo}

Let $Y$ be the alphabet $\{y_{n,\sigma}\ |\ (n,\sigma)\in\mathbb N_{>0}\times\Gamma\}$, indexed by $\mathbb N_{>0}\times\Gamma$. 
We define $\mathbb Q\langle Y\rangle$ to be the free algebra over this alphabet. We now 
construct linear maps fitting in the following diagram:  
$$
\xymatrix{
\mathbb Q\langle X\rangle\ar^{\pi_Y}@<2pt>[r]\ar_{\mathrm{corr}}@<-2pt>[r]\ar@/_2pc/_{(-)_*}[r] & \mathbb Q\langle Y\rangle
\ar_{\mathbf q}@(ur,ul)[]\ar^{\Delta_*}[r] & \mathbb Q\langle Y\rangle^{\otimes 2}}
$$ 

\begin{itemize}
\item the algebra $\mathbb Q\langle Y\rangle$ is viewed as a subalgebra of $\mathbb Q\langle X\rangle$ by the assignment 
$y_{n,\sigma}\mapsto x_0^{n-1}x_\sigma$ for $n\geq 0,\sigma\in\Gamma$. The map $\pi_Y$ is then determined by 
the conditions $\pi_Y(\mathbb Q\langle X\rangle x_0)=0$ and that $\mathbb Q\langle Y\rangle\hookrightarrow
\mathbb Q\langle X\rangle\stackrel{\pi_Y}{\to}\mathbb Q\langle Y\rangle$ is the identity (see \cite{Rac}, \S \S2.2.5 and 4.2.1); 

\item the map $\mathrm{corr}:\mathbb Q\langle X\rangle\to\mathbb Q\langle Y\rangle$ is defined to be the linear map taking 
the word $x_0^{n-1}x_1$ to ${{(-1)^{n-1}}\over n}(y_{1,1})^n$ if $n\geq 1$, and any other word to 0; in particular, it takes 1 to 0
(see \cite{Rac}, \S 3.1.1); 

\item the map $\mathbf q:\mathbb Q\langle Y\rangle\to\mathbb Q\langle Y\rangle$ is defined to be the linear map taking 1 to 1 
and by 
$$
\mathbf q(y_{s_1,\sigma_1}\cdots y_{s_r,\sigma_r})=y_{s_1,\sigma_1}y_{s_2,\sigma_2\sigma_1^{-1}}\cdots y_{s_r,\sigma_r
\sigma_{r-1}^{-1}}
$$
for any $r\geq 1$, $s_1,\ldots,s_r$ in $\mathbb N_{>0}$, and $\sigma_1,\ldots,\sigma_r$ in $\Gamma$
(see \cite{Rac}, \S 2.2.7). 

\item the map $(-)_*:\mathbb Q\langle X\rangle\to\mathbb Q\langle Y\rangle$ is given by $\psi\mapsto\psi_*$, where 
\begin{equation}\label{def:add:psi:star}
\psi_*:=\mathbf q(\pi_Y(\psi))+\mathrm{corr}(\psi); 
\end{equation}
one has therefore $(-)_*=\mathbf q\circ\pi_Y+\mathrm{corr}$
(see \cite{Rac}, (3.3.1.2)); 

\item the map $\Delta_*:\mathbb Q\langle Y\rangle\to\mathbb Q\langle Y\rangle^{\otimes 2}$ is the unique algebra morphism 
such that for any $n\geq 1$ and $\sigma\in\Gamma$ 
$$
\Delta_*(y_{n,\sigma})
=y_{n,\sigma}\otimes 1+1\otimes y_{n,\sigma}+\sum_{\stackrel{n'+n''=n|n',n''>0,}
{\stackrel{\sigma',\sigma''\in\Gamma|\sigma'\sigma''=\sigma}{}}}y_{n',\sigma'}\otimes y_{n'',\sigma''}. 
$$ 
The map $\Delta_*$ is called the  {\it harmonic coproduct}  (see \cite{Rac}, \S 2.3.1). 
\end{itemize}

\subsection{Modules over the Lie algebra $(\mathbb{Q}\langle X\rangle_0,\langle,\rangle)$}\label{sect:1:2:new}

The space $\mathbb Q\langle X\rangle$ is a left module over the Lie algebra $(\mathbb Q\langle X\rangle_0,\langle,\rangle)$
for the action
$$
\mathbb Q\langle X\rangle_0\to\mathrm{End}(\mathbb Q\langle X\rangle), \quad \psi\mapsto s_\psi
$$ 
(see \cite{Rac}, (3.1.9.2)). 

For $\varphi\in\mathbb Q\langle X\rangle_0$, there exists a unique linear endomorphism $s_\varphi^Y$ of 
$\mathbb Q\langle Y\rangle$ such that the following diagram commutes
\begin{equation}\label{diag:s:Y:phi}
\xymatrix{
\mathbb Q\langle X\rangle \ar^{s_\varphi}[r]\ar_{\mathbf q\circ\pi_Y}[d]& \mathbb Q\langle X\rangle\ar^{\mathbf q\circ\pi_Y}[d]\\ 
\mathbb Q\langle Y\rangle\ar_{s^Y_\varphi}[r]&\mathbb Q\langle Y\rangle}
\end{equation}
(see \cite{Rac}, \S 4.1.1). One shows: 
\begin{lemma}\label{lemma:tmqteqpts}
The map $\mathbb Q\langle X\rangle_0\to\mathrm{End}(\mathbb Q\langle Y\rangle)$, $\psi\mapsto s^Y_\psi$ defines 
a module structure on $\mathbb Q\langle Y\rangle$ over the Lie algebra $(\mathbb Q\langle X\rangle_0,\langle,\rangle)$. 
The map $\mathbb Q\langle X\rangle\stackrel{\mathbf q\circ\pi_Y}{\to}\mathbb Q\langle Y\rangle$ is then a module morphism
over this Lie algebra. 
\end{lemma}

\subsection{A Lie algebra homomorphism $\theta:(\mathfrak{Lib}(X),\langle,\rangle)
\to(\mathbb Q\langle X\rangle_0,\langle,\rangle)$}\label{sect:fctalah}

The following statement is a consequence of the fact that the elements $x_0^n,x_1^n$
are central in $(\mathbb Q\langle X\rangle_0,\langle,\rangle)$: 

\begin{lemma}\label{lemma:lfagbtcoc}
Let $(f_n)_{n>0}$ and $(g_n)_{n>0}$ be two collections of characters 
$f_n,g_n:(\mathfrak{Lib}(X),\langle,\rangle)\to\mathbb Q$. The map 
$$
\theta_{(f_n)_{n},(g_n)_{n}}:(\mathfrak{Lib}(X),\langle,\rangle)\to
(\mathbb Q\langle X\rangle_0,\langle,\rangle),\quad \psi\mapsto 
\psi+\sum_{n>0} f_n(\psi)x_1^n+\sum_{n>0} g_n(\psi)x_0^n
$$
is a Lie algebra homomorphism.
\end{lemma}

\begin{lemma}\label{lemma:tmfglqdb}
Let $n\geq 1$. Define\footnote{For $f$ an element of $\mathbb Q\langle X\rangle$ and $w$ a word in $X$, 
we denote by $(f|w)$ the coefficient of $w$ in $f$, so $f=\sum_w (f|w)w$. Viewing $\LL(X)$ as a subspace of $\mathbb Q\langle X\rangle$, 
we use the same notation for $f\in\LL(X)$.}  
maps $f_n^0,g_n^0:\mathfrak{Lib}(X)\to\mathbb Q$ by 
$f_n^0:\psi\mapsto(\psi|x_0^{n-1}x_1)$ for $n>0$, $g_1^0:\psi\mapsto (\psi|x_0)$, and $g_n^0:\psi\mapsto 0$ if $n>1$. 
These maps are Lie algebra characters $f_n^0,g_n^0:(\mathfrak{Lib}(X),\langle,\rangle)\to\mathbb Q$. 
\end{lemma}

{\em Proof.} According to Lemma \ref{lemma:taiozgla}, there is a decomposition 
$$
(\LL(X),\langle,\rangle)\simeq\mathfrak{outder}^*(\LL(X))^\Gamma\oplus(\mathbb Q x_0\oplus\mathbb Q x_1),  
$$ 
where the first summand is supported in bidegrees $(n,m)$ with $n,m>0$, and the second summand 
is supported in bidegrees $(1,0)$ and $(0,1)$. We will denote by $\mathrm{Supp}(X)$ the bidegree 
support of a bigraded vector space $X$. 

Let $n>1$. As $f_n^0$ is supported in bidegrees $(n-1,1)$ with $n-1>0$, this linear form factors
 through the summand $\mathfrak{outder}^*(\LL(X))^\Gamma$, and since 
$\mathrm{Supp}(\mathfrak{outder}^*(\LL(X))^\Gamma)\subset\mathbb Z_{>0}^2$,
one has $(n-1,1)\notin \mathrm{Supp}(\mathfrak{outder}^*(\LL(X))^\Gamma)+ 
\mathrm{Supp}(\mathfrak{outder}^*(\LL(X))^\Gamma)$, therefore $f_n^0$ is a character of the Lie 
algebra $\mathfrak{outder}^*(\LL(X))^\Gamma$, and therefore is a character of 
$(\LL(X),\langle,\rangle)$. 
By bidegree reasons, the linear forms $f_1^0$ and $g_1^0$ factor through the abelian Lie algebra 
$\mathbb Q x_0\oplus\mathbb Q x_1$, on which they are characters. These linear forms are 
therefore also characters of $(\LL(X),\langle,\rangle)$. \hfill\qed\medskip 

Combining Lemma \ref{lemma:lfagbtcoc} and Lemma \ref{lemma:tmfglqdb}, and
putting $\theta:=\theta_{(f_n^0)_n,(g_n^0)_n}$, we get: 
\begin{prop}\label{tmtltq}
The map $\theta:(\mathfrak{Lib}(X),\langle,\rangle)\to(\mathbb Q\langle X\rangle_0,\langle,\rangle)$ given by 
$$
\theta(\psi):=\psi+\sum_{n\geq 1}{(-1)^{n-1}\over n}(\psi|x_0^{n-1}x_1)x_1^n-(\psi|x_0)x_0
$$
is a Lie algebra homomorphism. 
\end{prop}

\subsection{The stabilizer Lie algebra of $\Delta_*$}\label{sect:stabilizer}

Via pull-back through the morphism $\theta$, the morphism 
$$
\mathbf q\circ\pi_Y : \mathbb{Q}\langle X\rangle\to\mathbb{Q}\langle Y\rangle
$$
of modules over $({\mathbb Q}\langle X\rangle_0,\langle,\rangle)$ may be viewed as a morphism of modules
over $({\LL}(X),\langle,\rangle)$ (see Lemma \ref{lemma:tmqteqpts}). The action of ${\LL}(X)$ on $\mathbb{Q}\langle X\rangle$ is given by 
$\psi\mapsto s_{\theta(\psi)}$, while its action on $\mathbb{Q}\langle Y\rangle$ is given by 
$\psi\mapsto s_{\theta(\psi)}^Y$.  

The space $\mathrm{Hom}_{\mathbb Q}(\mathbb Q\langle Y\rangle,\mathbb Q\langle Y\rangle^{\otimes 2})$ of linear maps 
$\mathbb Q\langle Y\rangle\to\mathbb Q\langle Y\rangle^{\otimes 2}$ is then equipped with a 
module structure over $(\LL(X),\langle,\rangle)$. Namely, the action of $\psi\in\LL(X)$ on $f\in
\mathrm{Hom}_{\mathbb Q}(\mathbb Q\langle Y\rangle,\mathbb Q\langle Y\rangle^{\otimes 2})$ is given by
$$\psi\cdot f:=
(s^Y_{\theta(\psi)}\otimes\mathrm{id}+\mathrm{id}\otimes s^Y_{\theta(\psi)})\circ f-f\circ s^Y_{\theta(\psi)}.$$ 
The stabilizer Lie algebra $\mathfrak{stab}(\Delta_*)$ of $\Delta_*\in
\mathrm{Hom}_{\mathbb Q}(\mathbb Q\langle Y\rangle,\mathbb Q\langle Y\rangle^{\otimes 2})$ is then a Lie subalgebra of 
$(\LL(X),\langle,\rangle)$, given by  
$$
\mathfrak{stab}(\Delta_*)=\{\psi\in\LL(X)\ |\ (s^Y_{\theta(\psi)}\otimes\mathrm{id}+\mathrm{id}\otimes
s^Y_{\theta(\psi)})\circ\Delta_*=\Delta_*\circ s^Y_{\theta(\psi)}\}.
$$

\subsection{Relation between $\tilde{\mathbf p}$, $\mathrm{sec}$ and $\theta$}\label{sect:relations}

Let $\partial_0$ be the derivation of $\mathbb Q\langle X\rangle$ defined by $\partial_0(x_0)=1$, $\partial_0(x_\sigma)=0$
for $\sigma\in\Gamma$ (see \cite{Rac}, \S 4.2.1). Let $\widetilde{\mathrm{sec}}$ be the linear endomorphism of 
$\mathbb Q\langle X\rangle$ given by 
$$
\widetilde{\mathrm{sec}}(f):=\sum_{i\geq 0}{(-1)^i\over i!}\partial_0^i(f)x_0^i
$$
for $f\in\mathbb Q\langle X\rangle$. 
Recall that $\mathbb Q\langle X\rangle$ contains subspaces $\mathbb Q\langle Y\rangle$ and $\mathrm{Ker}(\partial_0)$. 
It is proved in \cite{Rac}, \S 4.2.1, that $\widetilde{\mathrm{sec}}(\mathbb Q\langle Y\rangle)\subset \mathrm{Ker}(\partial_0)$. 
We denote by 
$$
\mathrm{sec}:\mathbb Q\langle Y\rangle\to\mathrm{Ker}(\partial_0)
$$
the resulting linear map. The following statement was proved in \cite{Rac}, Prop. 4.2.2: 

\begin{lemma}\label{lemma:inverse}
The map $\mathrm{sec}$ and the restriction $(\pi_Y)_{|\mathrm{Ker}(\partial_0)}:\mathrm{Ker}(\partial_0)\to
\mathbb Q\langle Y\rangle$ are inverse linear isomorphisms. 
\end{lemma}

Let $\mathbf p$ be the linear endomorphism of $\mathbb Q\langle Y\rangle$ given by 
\begin{equation}\label{def:p:Y}
\mathbf p(y_{n_1,\sigma_1}\cdots y_{n_r,\sigma_r})=y_{n_1,\sigma_1}y_{n_2,\sigma_1\sigma_2}\cdots y_{n_r,\sigma_1\cdots
\sigma_r}, 
\end{equation}
for any $n_i$ in $\mathbb N_{>0}$ and $\sigma_i$ in $\Gamma$. Then $\mathbf p$ and $\mathbf q$ are linear automorphisms of $\mathbb Q\langle Y\rangle$, which are inverse to each other (see \cite{Rac}, \S 2.2.7) . 

For $\psi\in\mathrm{Ker}(\partial_0)$, one has $\mathrm{sec}\circ\mathbf p(\psi_*)=\mathrm{sec}\circ\mathbf p(\mathbf q\circ 
\pi_Y(\psi)+\mathrm{corr}(\psi))$ (identity in $\mathrm{Ker}(\partial_0)$). Since $\mathbf p\circ\mathbf q
=\mathrm{id}_{\mathbb Q\langle Y\rangle}$ and by Lemma \ref{lemma:inverse}, one has 
$\mathrm{sec}\circ\mathbf p(\mathbf q\circ\pi_Y(\psi))=\psi$. The restriction of $\mathbf p$ and of $\mathrm{sec}$ to 
$\mathbb Q[y_{1,1}]$ is the identity.  As $\mathrm{corr}(\psi)$ belongs to $\mathbb Q[y_{1,1}]$, it follows that
$\mathrm{sec}\circ\mathbf p(\mathrm{corr}(\psi))=\mathrm{corr}(\psi)$. All this implies that  
\begin{equation}\label{8bis}
\forall\psi\in\mathrm{Ker}(\partial_0), \quad \mathrm{sec}\circ\mathbf p(\psi_*)=\psi+\mathrm{corr}(\psi) 
\end{equation}
(identity in $\mathrm{Ker}(\partial_0)$). 
Define $\LL(X)\ominus\mathbb Qx_0$ as the sum of all the homogeneous components of $\LL(X)$ for the degree
in $\mathbb N^{\{0\}\cup\Gamma}$, of degree $\neq \epsilon_0$. Then $\LL(X)\ominus\mathbb Qx_0\subset 
\mathrm{Ker}(\partial_0)$. One then has 
$$
\forall\psi\in\LL(X)\ominus\mathbb Qx_0, \quad \mathrm{sec}\circ\mathbf p(\psi_*)=\psi+\mathrm{corr}(\psi)
$$
(identity in $\mathrm{Ker}(\partial_0)$).
For any $\psi\in\LL(X)\ominus\mathbb QX_0$, one has $\psi+\mathrm{corr}(\psi)=\theta(\psi)$ (identity in 
$\mathbb Q\langle X\rangle$), therefore 
$$
\forall\psi\in \LL(X)\ominus\mathbb QX_0, \quad \mathrm{sec}\circ\mathbf p(\psi_*)=\theta(\psi)
$$
(identity in $\mathbb Q\langle X\rangle$).
For $\psi=x_0$, one has $\mathrm{sec}\circ\mathbf p(\psi_*)=\mathrm{sec}(0)=0=\theta(\psi)$. It follows: 
\begin{equation}\label{secp}
\forall\psi\in\LL(X), \quad 
\mathrm{sec}\circ\mathbf p(\psi_*)=\theta(\psi)
\end{equation}
(identity in $\mathbb Q\langle X\rangle$). 

Define $\tilde{\mathbf p}$ to be the endomorphism of $\mathbb Q\langle X\rangle$ given by 
$$
\tilde{\mathbf p}(x_0^{n_1}x_{\sigma_1}x_0^{n_2}\cdots x_0^{n_r}x_{\sigma_r}x_0^{n_{r+1}})=
x_0^{n_1}x_{\sigma_1}x_0^{n_2}\cdots x_0^{n_r}x_{\sigma_1\cdots\sigma_r}x_0^{n_{r+1}}, 
$$
for any $n_i$ in $\mathbb N$ and $\sigma_i$ in $\Gamma$. This endomorphism of $\mathbb Q\langle X\rangle$
commutes both with $\partial_0$ and with the endomorphism of $\mathbb Q\langle X\rangle$ given by right multiplication 
by $x_0$. It follows that $\tilde{\mathbf p}$ commutes with the endomorphism $\widetilde{\mathrm{sec}}$ of $\mathbb Q\langle X\rangle$. 
As the restriction of $\tilde{\mathbf p}$ to $\mathbb Q\langle Y\rangle\subset\mathbb Q\langle X\rangle$ is $\mathbf p$, it
follows that $\widetilde{\mathrm{sec}}\circ\mathbf{p}=\tilde{\mathbf{p}}\circ\mathrm{sec}$. Combining this with (\ref{secp}), one obtains: 

\begin{lemma}\label{prop:sec:sigma}
For any $\psi\in\LL(X)$, one has $\tilde{\mathbf p}\circ\mathrm{sec}(\psi_*)
=\theta(\psi)$ (identity in $\mathbb Q\langle X\rangle$). 
\end{lemma}

\section{Comparison of Lie algebras}\label{sect:cola}

This section is devoted to the comparison of the Lie algebras $\mathfrak{dmr}_0$ and $\mathfrak{stab}(\Delta_*)$. 
In \S \ref{sect:1:2}, we recall the definition of Racinet's spaces $\mathfrak{dmr}_0$ and $\mathfrak{dmr}$; we then  
compare the Lie algebra $\mathfrak{stab}(\Delta_*)$ introduced in \S \ref{sect:stabilizer} with these spaces. In \S \ref{sect:large:degrees}, 
we use the fact that $\mathfrak{dmr}_0$ and $\mathfrak{dmr}$ coincide in certain degrees to compute $\mathfrak{stab}(\Delta_*)$ 
in these degrees. In \S\S \ref{sect:comp:deg:1} and \ref{sect:cosdid2}, we compute 
$\mathfrak{stab}(\Delta_*)$ in the remaining degrees. This leads us in \S \ref{sect:main:thm} to Thm.\ \ref{main:thm} computing 
$\mathfrak{stab}(\Delta_*)$ in terms of $\mathfrak{dmr}_0$. We then derive Cor.\ \ref{cor:LA} stating that $\mathfrak{dmr}_0$ is a 
Lie subalgebra of $(\LL(X),\langle,\rangle)$, thus recovering Racinet's result (Prop.\ 4.1 i) in \cite{Rac}).

\subsection{The double inclusion $\mathfrak{dmr}_0\subset\mathfrak{stab}(\Delta_*)\subset\mathfrak{dmr}\oplus\mathbb Q x_0\oplus
\mathbb Q x_1$}\label{sect:double:incl}\label{sect:1:2} 

\begin{defin} (see \cite{Rac}, \S 3.3.1)
$$
\mathfrak{dmr}:=\{\psi\in\LL(X)\ |\ (\psi|x_0)=(\psi|x_1)=0, \quad \Delta_*(\psi_*)=\psi_*\otimes 1+1\otimes\psi_*\}
$$
\end{defin} 

\begin{defin} (see \cite{Rac}, \S 3.3.8)
One defines $\mathfrak{dmr}_0$ to be the subspace of $\mathfrak{dmr}$ of all elements $\psi$ 
satisfying the relation
$$
(\psi_*|y_{n,\sigma})+(-1)^n(\psi_*|y_{n,\sigma^{-1}})=0
$$
for $n=1$ and any $\sigma\in\Gamma$ if $|\Gamma|\geq 3$, and for $(n,\sigma)=(2,1)$ if $|\Gamma|\leq 2$. 
\end{defin}

In \cite{Rac}, Prop.\ 4.3.1 (see also \cite{Fur}, Prop.\ A.5 when $\Gamma=1$), it is proved: 
\begin{prop}\label{prop:Fur}
If $\psi\in\mathfrak{dmr}_0$, then 
\begin{equation}\label{id:Delta:star}
(s^Y_{\tilde{\mathbf p}\circ\mathrm{sec}(\psi_*)}\otimes\mathrm{id}+\mathrm{id}\otimes 
s^Y_{\tilde{\mathbf p}\circ\mathrm{sec}(\psi_*)})\circ\Delta_*=\Delta_*\circ s^Y_{\tilde{\mathbf p}\circ\mathrm{sec}(\psi_*)}
\end{equation}
(equality of maps $\mathbb Q\langle Y\rangle\to\mathbb Q\langle Y\rangle^{\otimes 2}$). 
\end{prop} 
By Lemma \ref{prop:sec:sigma}, condition (\ref{id:Delta:star}) is equivalent to condition 
$$(s^Y_{\theta(\psi)}\otimes\mathrm{id}+\mathrm{id}\otimes 
s^Y_{\theta(\psi)})\circ\Delta_*=\Delta_*\circ s^Y_{\theta(\psi)}$$
defining $\mathfrak{stab}(\Delta_*)$. 
Combining this with Prop.\ \ref{prop:Fur}, we get the inclusion 
\begin{equation}\label{iguazu1}
\mathfrak{dmr}_0\subset\mathfrak{stab}(\Delta_*). 
\end{equation}

Assume now that $\psi$ belongs to $\mathfrak{stab}(\Delta_*)$. It satisfies the identity 
$$
(s^Y_{\theta(\psi)}\otimes\mathrm{id}+\mathrm{id}\otimes s^Y_{\theta(\psi)})\circ\Delta_*=\Delta_*\circ s^Y_{\theta(\psi)}.
$$ 
Applying this identity to $1\in\mathbb Q\langle Y\rangle$, one obtains 
\begin{equation}\label{id:sec:psi:*}
s^Y_{\theta(\psi)}(1)\otimes1+1\otimes s^Y_{\theta(\psi)}(1)=\Delta_*( s^Y_{\theta(\psi)}(1))
\end{equation}
(identity in $\mathbb Q\langle Y\rangle^{\otimes 2}$). 
Let $\varphi$ be an element of $\mathbb Q\langle X\rangle$. 
Applying diagram (\ref{diag:s:Y:phi}) to $1\in\mathbb Q\langle X\rangle$, one obtains 
$s_{\varphi}^Y(\mathbf q(\pi_Y(1)))=\mathbf q\circ\pi_Y(s_\varphi(1))$, which using $\pi_Y(1)=\mathbf q(1)=1$ 
and $s_\varphi(1)=\varphi$, implies: 
$$
\forall\varphi\in\mathbb Q\langle X\rangle, \quad s_\varphi^Y(1)=\mathbf q\circ\pi_Y(\varphi). 
$$
Applying this identity to $\varphi :=\theta(\psi)$, one obtains $s_{\theta(\psi)}^Y(1)
={\mathbf q}\circ\pi_Y(\theta(\psi))$. By Lemma \ref{prop:sec:sigma}, this is equal to 
${\mathbf q}\circ\pi_Y(\tilde{\mathbf p}\circ\mathrm{sec}(\psi_*))$. We have $\pi_Y\circ\tilde{\mathbf p}
=\mathbf p\circ\pi_Y$; moreover $\mathbf q\circ\mathbf p=\mathrm{id}_{\mathbb Q\langle Y\rangle}$, and Lemma \ref{lemma:inverse} 
implies that $\pi_Y\circ\mathrm{sec}=\mathrm{id}_{\mathbb Q\langle Y\rangle}$, so that 
$s_{\theta(\psi)}^Y(1)=\psi_*$. 

%
Identity (\ref{id:sec:psi:*}) then implies that 
$$\psi_*\otimes 1+1\otimes\psi_*=\Delta_*(\psi_*)$$ (identity in $\mathbb Q\langle Y\rangle^{\otimes 2}$). 
Therefore: 
\begin{equation}\label{iguazu2}
\mathfrak{stab}(\Delta_*)\subset \{\psi\in\LL(X)\ |\ \Delta_*(\psi_*)=\psi_*\otimes 1+1\otimes\psi_*\}.
\end{equation} 
(\ref{iguazu1}) and (\ref{iguazu2}) then imply
\begin{equation}\label{double:inclusion}
\mathfrak{dmr}_0\subset\mathfrak{stab}(\Delta_*)\subset\mathfrak{dmr}\oplus\mathbb Q x_0\oplus
\mathbb Q x_1.
\end{equation}

\subsection{Computation of $\mathfrak{stab}(\Delta_*)$ for large degrees}\label{sect:large:degrees}

The double inclusion (\ref{double:inclusion}) deals with vector subspaces of the space $\LL(X)$. This space is equipped with a grading for 
which the generators $x_0$ and $x_\sigma$ for $\sigma\in\Gamma$ have degree 1; this is the total degree of the bidegree from 
\S\ref{sect:odla}. The spaces from (\ref{double:inclusion}) are then graded subspaces
of $\LL(X)$. The space $\LL(X)$ is the direct sum of its components of degrees $\geq 1$ and the extremal spaces
$\mathfrak{dmr}_0$ and $\mathfrak{dmr}\oplus\mathbb Q x_0\oplus\mathbb Q x_1$ of the double inclusion (\ref{double:inclusion})
coincide in each degree $\geq 2$ if $|\Gamma|\geq 3$, and in each degree $\geq 3$ if $|\Gamma|=1,2$ (here $|\Gamma|$ is the cardinality
of $\Gamma$). 

It follows: 

\begin{lemma}\label{lemma:large:degrees}
\begin{itemize}
\item If $|\Gamma|\geq 3$, then $\mathfrak{stab}(\Delta_*)[n]=\mathfrak{dmr}_0[n]$ any $n\geq 2$. 

\item  If $|\Gamma|=1,2$, then $\mathfrak{stab}(\Delta_*)[n]=\mathfrak{dmr}_0[n]$ any $n\geq 3$. 
\end{itemize}
\end{lemma} 
Here the notation $[n]$ means the component of degree $n$. 

\subsection{Computation of $\mathfrak{stab}(\Delta_*)$ in degree $1$}\label{sect:comp:deg:1}

The degree 1 component of $\mathfrak{dmr}\oplus\mathbb Q x_0\oplus\mathbb Q x_1$ coincides with that of $\LL(X)$, which is 
$$
\oplus_{g\in\Gamma\cup\{0\}}\mathbb Q x_g.
$$
On the other hand, according to \S\ref{sect:1:2}, the degree 1 component of
$\mathfrak{dmr}_0$ is 
$$
\mathrm{Span}\{x_g+x_{g^{-1}}\ | \  g\in\Gamma\setminus\{1\}\}.
$$ 

Recall that $\mathbb Q\langle Y\rangle$ may be viewed as the subspace $\mathbb Q\oplus\oplus_{\alpha\in\Gamma}
\mathbb Q\langle X\rangle x_\alpha$ of $\mathbb Q\langle X\rangle$ (see \S\ref{sect:hccpo}). 
If now $\psi$ belongs to $\mathbb Q\langle Y\rangle$ and $\alpha$ belongs to $\Gamma$, the element 
$d_{\psi}(x_\alpha)=[x_\alpha,t_\alpha(\psi)]$ also belongs to $\mathbb Q\langle Y\rangle$.
As $d_{\psi}$ is a derivation, it follows that $d_{\psi}$ preserves $\mathbb Q\langle Y\rangle$. 
We denote by 
$$
\tilde d_{\psi}:\mathbb Q\langle Y\rangle\to\mathbb Q\langle Y\rangle
$$
the resulting linear map. On the other hand, the endomorphism of $\mathbb Q\langle X\rangle$ given by $v\mapsto \psi v$
also preserves the subspace $\mathbb Q\langle Y\rangle$. As the endomorphism $s_{\psi}$ of $\mathbb Q\langle X\rangle$
is the sum of $d_{\psi}$ and of $(v\mapsto\psi v)$, it follows that $s_{\psi}$ preserves $\mathbb Q\langle Y\rangle$.
We denote by 
$$
\tilde s_{\psi}:\mathbb Q\langle Y\rangle\to\mathbb Q\langle Y\rangle
$$
the resulting linear map. The purpose of this section is to compute the degree 1 component of the space 
$\mathfrak{stab}(\Delta_*)$, which lies between the extremal spaces of (\ref{double:inclusion}). 

\subsubsection{$s^Y_{\psi}$ when $\psi$ lies in $\mathbb Q\langle Y\rangle$}

Combining the diagram describing the compatibility of $s_{\psi}$ and $\tilde s_{\psi}$ with the diagram 
(\ref{diag:s:Y:phi}), we obtain the diagram 
$$
\xymatrix{\mathbb Q\langle Y\rangle\ar@{^(->}[r]\ar_{\tilde s_{\psi}}[d]& 
\mathbb Q\langle X\rangle\ar_{s_{\psi}}[d]\ar^{\mathbf q\circ\pi_Y}[r] 
& \mathbb Q\langle Y\rangle\ar^{s_{\psi}^Y}[d] \\ \mathbb Q\langle Y\rangle\ar@{^(->}[r]
& \mathbb Q\langle Y\rangle\ar^{\mathbf q\circ\pi_Y}[r]& \mathbb Q\langle Y\rangle}
$$
Since the composed map $\mathbb Q\langle Y\rangle\hookrightarrow\mathbb Q\langle X\rangle\stackrel{\pi_Y}{\to}
\mathbb Q\langle Y\rangle$ is the identity of $\mathbb Q\langle Y\rangle$, the composition of the horizontal maps
of this diagram equals $\mathbf q$, so that 
\begin{equation}\label{rel:tildes:sY}
\forall\psi\in\mathbb Q\langle Y\rangle, \quad s_{\psi}^Y=\mathbf q\circ\tilde s_{\psi}\circ\mathbf q^{-1}.
\end{equation}

\subsubsection{Condition for $s^Y_\psi$ to ba a coderivation}
For $g\in\Gamma$, one has 
$$
s_{x_g}(1)=x_g\quad \mathrm{and} \quad \forall\alpha\in\Gamma, \quad 
s_{x_g}(x_\alpha)=x_gx_\alpha+[x_\alpha,x_{\alpha g}]. 
$$
As $x_g$ lies in $\mathbb Q\langle Y\rangle$ (it identifies with the element $y_{1,g}$), the endomorphism  
$s_{x_g}$ of $\mathbb Q\langle X\rangle$ restricts to an endomorphism $\tilde s_{x_g}$ of $\mathbb Q\langle Y\rangle$, 
which satisfies 
$$
\tilde s_{x_g}(1)=y_{1,g}\quad \mathrm{and} \quad \forall\alpha\in\Gamma, \quad 
\tilde s_{x_g}(y_{1,\alpha})=y_{1,g}y_{1,\alpha}+[y_{1,\alpha},y_{1,\alpha g}]. 
$$
Equation (\ref{rel:tildes:sY}) then implies: 
$$
s_{x_g}^Y(1)=y_{1,g}\quad\mathrm{and}\quad
\forall\alpha\in\Gamma, \quad 
s_{x_g}^Y(y_{1,\alpha})=y_{1,g}y_{1,\alpha g^{-1}}+y_{1,\alpha}y_{1,g}-y_{1,\alpha g}y_{1,g^{-1}}. 
$$
Let $\alpha\in\Gamma$. As $y_{1,\alpha}$ is primitive for $\Delta_*$, one has 
$$
(s_{x_g}^Y\otimes\mathrm{id}+\mathrm{id}\otimes s_{x_g}^Y)\circ\Delta_*(y_{1,\alpha})
=s_{x_g}^Y(y_{1,\alpha})\otimes 1+1\otimes s_{x_g}^Y(y_{1,\alpha})+\mathrm{sym}(y_{1,\alpha}\otimes y_{1,g}),   
$$
where $\mathrm{sym}(a\otimes b):=a\otimes b+b\otimes a$. On the other hand, 
\begin{align*}
\Delta_*(s_{x_g}^Y(y_{1,\alpha}))=& s_{x_g}^Y(y_{1,\alpha})\otimes 1 
+1\otimes s_{x_g}^Y(y_{1,\alpha}) \\
& \quad +\mathrm{sym}(y_{1,\alpha}\otimes y_{1,g})+\mathrm{sym}(y_{1,g}\otimes y_{1,\alpha g^{-1}})
-\mathrm{sym}(y_{1,g^{-1}}\otimes y_{1,\alpha g}). 
\end{align*}
It follows that 
\begin{equation}\label{identity:y:1:alpha}
\big(\Delta_*\circ s_{x_g}^Y-(s_{x_g}^Y\otimes\mathrm{id}+\mathrm{id}\otimes s_{x_g}^Y)\circ\Delta_*\big)(y_{1,\alpha})
=\mathrm{sym}(y_{1,g}\otimes y_{1,\alpha g^{-1}})-\mathrm{sym}(y_{1,g^{-1}}\otimes y_{1,\alpha g}). 
\end{equation}

Let now $\psi$ be an element of $\mathrm{Span}\{x_g|g\in\Gamma\}$ and assume that $s^Y_{\psi}$
is a coderivation for $\Delta_*$, namely that $\Delta_*\circ s_{\psi}^Y=(s_{\psi}^Y\otimes\mathrm{id}+\mathrm{id}\otimes 
s_{\psi}^Y)\circ\Delta_*$. Set $\psi=\sum_{g\in\Gamma}c_g x_g$. Using (\ref{identity:y:1:alpha}), the equation 
$\big(\Delta_*\circ s_{\psi}^Y-(s_{\psi}^Y\otimes\mathrm{id}+\mathrm{id}\otimes 
s_{\psi}^Y)\circ\Delta_*\big)(y_{1,\alpha})=0$ expresses as follows
$$
\forall \alpha\in\Gamma, \quad \sum_{g\in\Gamma}c_g(x_g\otimes\alpha x_{g^{-1}}-x_{g^{-1}}\otimes x_{\alpha g}
+x_{\alpha g^{-1}}\otimes x_g
-x_{\alpha g}\otimes x_{g^{-1}})=0
$$
(equality in $\mathbb Q\Gamma^{\otimes 2}$). Set $A:=\sum_{g\in\Gamma}c_g(x_g\otimes x_{g^{-1}}-x_{g^{-1}}\otimes x_g)$. Multiplying 
the above identity by $1\otimes x_{\alpha^{-1}}$, it reexpresses as 
\begin{equation}\label{cond:A}
\forall\alpha\in\Gamma,\quad A=(x_\alpha\otimes x_{\alpha^{-1}})A. 
\end{equation}
For any $g\in\Gamma$, set $a_g:=c_g-c_{g^{-1}}$. Then $A=\sum_{g\in\Gamma}a_g(x_g\otimes x_{g^{-1}})$. As $\Gamma$
is abelian, (\ref{cond:A}) translates into the condition $a_g=a_{\alpha^{-1}g}$ for any $(g,\alpha)\in\Gamma^2$. This means
that the map $g\mapsto a_g$ is constant. On the other hand, $a_1=0$, therefore $a_g=0$ for any $g\in\Gamma$. This implies that 
$c_g=c_{g^{-1}}$ for any $g\in\Gamma$. We have proved:

\begin{lemma}\label{lemma:coder}
If $\psi$ is an element of $\mathrm{Span}\{x_g | g\in\Gamma\}$ such that $s^Y_{\psi}$
is a coderivation for $\Delta_*$, then it belongs to $\mathrm{Span}\{x_g+x_{g^{-1}} | g\in\Gamma\}$.
\end{lemma} 

The restriction to $X$ of the map $\theta:\LL(X)\to\mathbb Q\langle X\rangle_0$ is given by 
\begin{equation}\label{eqs:sigma}
x_0\mapsto 0,\quad x_1\mapsto 2x_1, \quad x_g\mapsto x_g\quad\mathrm{for}\quad g\in\Gamma\setminus\{1\}.
\end{equation}
The degree 1 part of 
$\mathfrak{stab}(\Delta_*)$ is the set of all $\underline\psi\in\mathrm{Span}\{x_g|g\in\Gamma\}$, such that
$s^Y_{\theta(\underline\psi)}$ is a coderivation. Combining Lemma \ref{lemma:coder} with (\ref{eqs:sigma}), 
we obtain: 

\begin{lemma}\label{lemma:stab:1}
The degree $1$ part of $\mathfrak{stab}(\Delta_*)$ is contained in $\mathbb Q x_0\oplus\mathrm{Span}\{x_g+x_{g^{-1}}|
g\in\Gamma\}$.  
\end{lemma}
 
On the other hand, it follows from the definition of $\mathfrak{dmr}_0$ that its degree 1 part is equal to 
$$
\mathrm{Span}\{x_g+x_{g^{-1}}|g\in\Gamma\setminus\{1\}\}
$$
(see \S\ref{sect:1:2}). Combining this equality, the first inclusion of (\ref{double:inclusion}) and 
Lemma \ref{lemma:stab:1}, we obtain: 
\begin{lemma}\label{lemma:comp:stab:deg:1}
Let $\Gamma$ be arbitrary. The components of degree $1$ of $\mathfrak{stab}(\Delta_*)$ and $\mathfrak{dmr}_0$
are related by 
$$
\mathfrak{stab}(\Delta_*)[1]=\mathfrak{dmr}_0[1]\oplus\mathbb Q x_0\oplus\mathbb Q x_1. 
$$
\end{lemma} 

\subsection{Computation of $\mathfrak{stab}(\Delta_*)$ in degree $2$}\label{sect:cosdid2}

\subsubsection{The case $|\Gamma|=1$}\label{sect:Gamma=1:n=2}

Specializing the inclusion $\mathfrak{stab}(\Delta_*)\subset \mathfrak{dmr}\oplus\mathbb Q x_0\oplus\mathbb Q x_1\subset\LL(X)$ 
in degree $2$, we obtain the inclusion 
$$
\mathfrak{stab}(\Delta_*)[2]\subset \mathfrak{dmr}[2]\subset\LL(X)[2]
$$
as $x_0,x_1$ are of degree 1. 

One has $\LL(X)[2]=\mathbb Q\cdot[x_0,x_1]$. If $\psi:=[x_0,x_1]$, then $\psi_*=y_2-{1\over 2}y_1^2$ (we drop the elements of 
$\Gamma$ from the notation of generators of $\mathbb Q\langle Y\rangle$ as $\Gamma=1$), which is $\Delta_*$-primitive. 
Therefore $\mathfrak{dmr}[2]=\mathbb Q\cdot [x_0,x_1]$. 

One computes $\sigma([x_0,x_1])=[x_0,x_1]-{1\over 2}x_1^2$. The image of $x_1$ by $s_{\sigma([x_0,x_1])}$ is 
$$
s_{\sigma([x_0,x_1])}(x_1)=\sigma([x_0,x_1]) x_1+[x_1,\sigma([x_0,x_1])]=x_1\sigma([x_0,x_1])=x_1[x_0,x_1]-{1\over 2}x_1^3.
$$ The image of this element under $\mathbf q\circ\pi_Y$ is 
$y_1y_2-{1\over 2}y_1^3$. 
So 
$$
\mathbf q\circ\pi_Y\circ s_{\sigma([x_0,x_1])}(x_1)=y_1y_2-{1\over 2}y_1^3.
$$
On the other hand, the image of $x_1$ by $\mathbf q\circ\pi_Y$ is $y_1$. It then follows from diagram 
(\ref{diag:s:Y:phi}) that 
$$
s^Y_{\sigma([x_0,x_1])}(y_1)=y_1y_2-{1\over 2}y_1^3. 
$$
Whereas $y_1$ is primitive for $\Delta_*$, the element $y_1y_2-{1\over 2}y_1^3$ is not, as its image by $\Delta_*$ contains
for example the term $y_1\otimes y_2$. It follows that $s^Y_{\sigma([x_0,x_1])}$ is not a coderivation for $\Delta_*$, so that 
$[x_0,x_1]\notin\mathfrak{stab}(\Delta_*)[2]$. Therefore 
\begin{lemma}\label{Gamma=1:n=2=stab:dmr0}
When $|\Gamma|=1$, one has 
$$
\mathfrak{stab}(\Delta_*)[2]=0=\mathfrak{dmr}_0[2]. 
$$
\end{lemma}

\subsubsection{The case $|\Gamma|=2$}\label{sect:Gamma=1:n=2}\label{sect:Gamma=2:n=2}

Specializing the inclusion $\mathfrak{stab}(\Delta_*)\subset \mathfrak{dmr}\oplus\mathbb Q x_0\oplus\mathbb Q x_1\subset\LL(X)$ 
in degree $2$, we again obtain the inclusion 
$$
\mathfrak{stab}(\Delta_*)[2]\subset \mathfrak{dmr}[2]\subset\LL(X)[2]. 
$$

We set $\Gamma=\{+,-\}$ so that $x_1$ can be denoted $x_+$. One has 
$$
\LL(X)[2]=\mathbb Q[x_0,x_+]\oplus \mathbb Q[x_0,x_-]\oplus \mathbb Q[x_+,x_-]. 
$$
Set $$
\psi_1:=[x_0,x_+], \quad \psi_2:=[x_0,x_-],\quad \psi_3:=[x_+,x_-].
$$
Then $(\psi_1)_*=y_{2,+}-{1\over 2}(y_{1,+})^2$, 
$(\psi_2)_*=y_{2,-}$, $(\psi_3)_*=y_{1,+}y_{1,-}-(y_{1,-})^2$. Then 
\begin{align*}
\Delta_\ast((\psi_1)_*)&=(\psi_1)_*\otimes 1+1\otimes(\psi_1)_*+y_{1,-}\otimes y_{1,-}, \\ 
\Delta_\ast((\psi_2)_*)&=(\psi_2)_*\otimes 1+1\otimes(\psi_2)_*+y_{1,+}\otimes y_{1,-}+y_{1,-}\otimes y_{1,+},\\ 
\Delta_\ast((\psi_3)_*)&=(\psi_3)_*\otimes 1+1\otimes(\psi_3)_*
+y_{1,+}\otimes y_{1,-}+y_{1,-}\otimes y_{1,+}-2y_{1,-}\otimes y_{1,-}.
\end{align*}
A linear combination of $(\psi_1)_*,(\psi_2)_*,(\psi_3)_*$ is therefore primitive for $\Delta_*$ iff it is a linear combination of 
$2(\psi_1)_*-(\psi_2)_*+(\psi_3)_*$. It follows that 
$$
\mathfrak{dmr}[2]=\mathbb Q\cdot \psi, \quad\mathrm{where}\quad \psi:=2\psi_1-\psi_2+\psi_3. 
$$
One computes $\sigma(\psi)=\psi+{(-1)^1\over 2}(\psi|x_0x_+)x_+^2$, therefore 
$$
\sigma(\psi)=\psi-x_+^2.
$$ 
The image of $x_+$ by $s_{\sigma(\psi)}$ is $\sigma(\psi)x_++d_{\sigma(\psi)}(x_+)=x_+\sigma(\psi)$ so 
$$
s_{\sigma(\psi)}(x_+)=2(x_+x_0x_+ - x_+x_+x_0)-(x_+x_0x_- - x_+x_-x_0)+(x_+x_+x_- - x_+x_-x_+). 
$$
The image of this element by $\pi_Y$ is 
$$
2y_{1,+}y_{2,+}-y_{1,+}y_{2,-} +(y_{1,+})^2y_{1,-}-y_{1,+}y_{1,-}y_{1,+}.
$$
The image of this element by $\mathbf q$ is 
$$
2y_{1,+}y_{2,+}-y_{1,+}y_{2,-} +(y_{1,+})^2y_{1,-}-y_{1,+}(y_{1,-})^2.
$$
Since the image of $x_+$ by $\mathbf q\circ\pi_Y$ is $y_{1,+}$, (\ref{diag:s:Y:phi}) implies 
$$
s_{\sigma(\psi)}^Y(y_{1,+})=
2y_{1,+}y_{2,+}-y_{1,+}y_{2,-} +(y_{1,+})^2y_{1,-}-y_{1,+}(y_{1,-})^2.
$$
Whereas $y_{1,+}$ is primitive for $\Delta_*$, its image by $s_{\sigma(\psi)}^Y$ is not: the image by $\Delta_*$ of 
$s_{\sigma(\psi)}^Y(y_{1,+})$ contains for example the word $y_{1,+}\otimes y_{2,+}$. It follows that $s_{\sigma(\psi)}^Y$ is 
not a coderivation for $\Delta_*$, so that $\psi\notin\mathfrak{stab}(\Delta_*)[2]$. Therefore:
\begin{lemma}\label{Gamma=2:n=2=stab:dmr0}
When $|\Gamma|=2$, one has
$$
\mathfrak{stab}(\Delta_*)[2]=0=\mathfrak{dmr}_0[2].
$$  
\end{lemma}

\subsection{The equality $\mathfrak{stab}(\Delta_*)=\mathfrak{dmr}_0\oplus\mathbb Q x_0\oplus
\mathbb Q x_1$}\label{sect:main:thm}

Combining Lemma \ref{lemma:large:degrees}, Lemma \ref{lemma:comp:stab:deg:1}, Lemma \ref{Gamma=1:n=2=stab:dmr0}, and 
Lemma \ref{Gamma=2:n=2=stab:dmr0}, we obtain: 

\begin{thm}\label{main:thm}
Let $\Gamma$ be arbitrary. The Lie algebras $\mathfrak{stab}(\Delta_*)$ and $\mathfrak{dmr}_0\oplus
\mathbb Q x_0\oplus\mathbb Q x_1$ are equal:  
$$
\mathfrak{stab}(\Delta_*)=\mathfrak{dmr}_0\oplus\mathbb Q x_0\oplus\mathbb Q x_1. 
$$
\end{thm} 
Note that $\mathbb Q x_0\oplus\mathbb Q x_1$ is contained in the center of $(\LL(X),\langle,\rangle)$. 
The subspace $\mathfrak{stab}(\Delta_*)$ of $(\LL(X),\langle,\rangle)$ is obviously a Lie subalgebra. 
It is graded and lives in degrees $\geq 1$. Therefore projection to any quotient of its degree 1 component is a Lie 
algebra morphism to an abelian Lie algebra, and the kernel of such a projection is an ideal in $\mathfrak{stab}(\Delta_*)$, 
and therefore a Lie subalgebra of $(\LL(X),\langle,\rangle)$.  Applying this to the quotient of $\mathfrak{stab}(\Delta_*)[1]$
by its subspace $\mathbb Q x_0\oplus\mathbb Q x_1$, we obtain: 

\begin{cor}\label{cor:LA}
$\mathfrak{dmr}_0$ is a Lie subalgebra of 
$(\LL(X),\langle,\rangle)$. 
\end{cor} 
As mentioned in the introduction, this result was first proved in \cite{Rac}. 

\begin{remark}
Let $\Gamma,\Gamma'$ be finite commutative groups. Let us specify the dependence of 
$X$ on the group $\Gamma$ by denoting it $X_\Gamma$. 
To a group morphism $\varphi:\Gamma\to\Gamma'$, one associates the algebra morphisms 
$$
\varphi_* : \mathbb Q\langle X_\Gamma\rangle\to \mathbb Q\langle X_{\Gamma'}\rangle
\quad \mathrm{and}\quad 
\varphi^* : \mathbb Q\langle X_{\Gamma'}\rangle\to \mathbb Q\langle X_\Gamma\rangle, 
$$
defined by 
$$
\varphi_*:x_0\mapsto d\cdot x_0, \quad x_\gamma\mapsto x_{\varphi(\gamma)} \quad\mathrm{for}\quad \gamma\in\Gamma, 
$$
where $d:=|\mathrm{Ker}\varphi|$, and 
$$
\varphi^*:x_0\mapsto x_0, \quad x_{\gamma'}\mapsto\sum_{\gamma\in\varphi^{-1}(\gamma')} x_\gamma\quad\mathrm{for}\quad \gamma'\in\Gamma'.  
$$

The maps $\varphi_*:\mathbb Q\langle X_\Gamma\rangle\to\mathbb Q\langle X_{\Gamma'}\rangle$ and 
$\varphi^*:\mathbb Q\langle X_{\Gamma'}\rangle\to\mathbb Q\langle X_\Gamma\rangle$ are Lie algebra 
morphisms, the spaces $\mathbb Q\langle X_\Gamma\rangle$ and $\mathbb Q\langle X_{\Gamma'}\rangle$ being equipped with the 
Ihara brackets. As the maps $\varphi_*$ and $\varphi^*$ restrict to linear maps $\LL(X_\Gamma)\rightleftarrows\LL(X_{\Gamma'})$, the maps $\varphi_*$ and $\varphi^*$ define Lie algebra morphisms between the Ihara Lie algebras 
$(\LL(X_\Gamma),\langle,\rangle)$ and $(\LL(X_{\Gamma'}),\langle,\rangle)$. 

If $\Gamma$ is a finite commutative group, then for any integer $d\geq 1$, the group $\Gamma^d$ is defined
as the image of the endomorphism of $\Gamma$ given by multiplication by $d$. There is an injective morphism 
$i_d:\Gamma^d\to\Gamma$ given by inclusion, and a surjective morphism $p^d:\Gamma\to\Gamma^d$ given by corestriction of 
the multiplication by $d$. To them are associated morphisms $i_d^*,p^d_*:\LL(X_\Gamma)\to\LL(X_{\Gamma^d})$ of Lie algebras, 
both sides being equipped with the Ihara bracket. Any linear form $\ell_d:\mathbb QX_\Gamma\to\mathbb Q$ extends to a linear 
form $\tilde\ell_d:\LL(X_\Gamma)\to\mathbb Q$ by the condition that $\tilde\ell_d$ vanishes on the components of degree $>1$. 
As the image of the Ihara Lie bracket is concentrated in degrees $>1$, the form $\tilde\ell_d$ is a character of the Ihara bracket on 
$\LL(X_\Gamma)$. Since on the other hand $x_0\in\LL(X_\Gamma)$ is a central element, the space 
$$
\LL(X_\Gamma)_d:=\{x\in\LL(X_\Gamma)\ |\ p_*^d(x)=i^*_d(x)+\tilde\ell_d(x)x_0\}
$$
is a Lie subalgebra of $(\LL(X_\Gamma),\langle,\rangle)$. Set 
$\ell_d(x):=\sum_{\sigma\in\Gamma|\sigma^{n/d}=1}(x|x_\sigma)$, and 
$$
\LL(X_\Gamma)_{\mathrm{dist}}:=\bigcap_{d||\Gamma|}\LL(X_\Gamma)_d. 
$$
Then $\LL(X_\Gamma)_{\mathrm{dist}}$ is a Lie subalgebra of $\LL(X_\Gamma)$. The Lie algebra $\mathfrak{dmrd}$ from \cite{Rac}, 
\S 3.3.1, is the intersection $\LL(X_\Gamma)_{\mathrm{dist}}\cap\mathfrak{dmr}_0$.  \hfill\qed\medskip 
\end{remark}


\section{Exponentiation of the Lie algebra homomorphism $\theta$}\label{sect:eotlaht}

Thm. \ref{main:thm} says that $\mathfrak{dmr}_0$ identifies with the stabilizer of an element in a module over 
$({\LL}(X),\langle,\rangle)$, obtained from a $\mathfrak{mt}$-module by the pull-back under $\theta:({\LL}(X),\langle,
\rangle)\to\mathfrak{mt}$. In order to prove the group version of this result, we construct in this section the group analogue 
of $\theta$. In \S\ref{sect:qgsatla}, we recall useful notions from the theory of $\mathbb Q$-group schemes. In \S\ref{sect:aqgsmtm}, 
we show how a collection of characters of the group scheme $(\Bbbk\mapsto(\mathrm{exp}(\hat{\LL}_\Bbbk(X)),\circledast))$ gives rise to a morphism of this group scheme to $\sf{MT}$. In \S\ref{sect:coaqgsmtet}, we construct an explicit collection of such characters, giving 
rise to a morphism $\Theta$ to $\sf MT$, whose associated Lie algebra morphism is $\theta$ (Prop. \ref{lemma:constr:Theta}). 

\subsection{$\mathbb Q$-group schemes and their Lie algebras}\label{sect:0:3}\label{sect:qgsatla}

\subsubsection{} 

An affine $\mathbb Q$-group scheme is a functor $\{$commutative $\mathbb Q$-algebras$\}\to\{$groups$\}$ which is 
representable by a commutative Hopf algebra. Such a functor ${\sf G}:\{$commutative $\mathbb Q$-algebras$\}\to\{$groups$\}$ being fixed, 
the kernel 
$\mathrm{Ker}({\sf G}(\mathbb Q[\epsilon]/(\epsilon^2))\to {\sf G}(\mathbb Q))$ is naturally equipped with a $\mathbb Q$-Lie algebra structure; 
this is the Lie algebra $(\mathrm{Lie}({\sf G}),[,])$ of ${\sf G}$. If ${\sf G},{\sf H}$ are two $\mathbb Q$-group schemes, then a morphism 
${\sf G}\to{\sf H}$
is a morphism from the Hopf algebra of ${\sf H}$ to that of ${\sf G}$. It gives rise to a natural transformation between the functors 
$\{$commutative $\mathbb Q$-algebras$\}\to\{$groups$\}$ attached to ${\sf G},{\sf H}$.

\subsubsection{} \label{rtaqla}

Recall that a $\mathbb Q$-Lie algebra is called nilpotent if its lower central series stabilizes to 0. 
A $\mathbb Q$-Lie algebra is called pronilpotent if it is complete and separated for the topology defined by the members of the 
lower central series. Such a Lie algebra is then a projective limit of nilpotent Lie algebras, which are its quotients by the members of 
the lower central series. 

If $\mathfrak g$ is a pronilpotent $\mathbb Q$-Lie algebra and $\mathbf k$ is a commutative $\mathbb Q$-algebra, define $\mathfrak g\hat\otimes\mathbf k$ to be the inverse limit $\lim_{\leftarrow}(\mathfrak g/\mathfrak g^{(n)})\otimes\mathbf k$,  where $\mathfrak g^{(n)}:=
[\mathfrak g,[\mathfrak g,\ldots,\mathfrak g]]$ ($n$ arguments). This is a complete separated topological Lie algebra over $\mathbf k$. 

For $x,y\in\mathfrak g\hat\otimes\mathbf k$, the Campbell-Baker-Hausdorff (CBH) series $\mathrm{cbh}(x,y)=\mathrm{log}(e^xe^y)$ 
is a well-defined element of 
$\mathfrak g\hat\otimes\mathbf k$. When equipped with the product map $\mathrm{cbh}:(\mathfrak g\hat\otimes\mathbf k)^2\to\mathfrak g\hat\otimes\mathbf k$, $(x,y)\mapsto \mathrm{cbh}(x,y)$, the set $\mathfrak g\hat\otimes\mathbf k$ is a group. The assignment
$\{$commutative $\mathbb Q$-algebras$\}\to\{$groups$\}$, $\mathbf k\mapsto(\mathfrak g\hat\otimes\mathbf k,\mathrm{cbh})=:
{\sf G}_{\mathfrak n}(\Bbbk)$ is a $\mathbb Q$-group scheme, denoted ${\sf G}_{\mathfrak n}$ and 
represented by the Hopf algebra given by the direct limit $\limm_{\rightarrow}U(\mathfrak g/\mathfrak g^{(n)})^\circ$, where $U$
is the universal enveloping algebra functor and $^\circ$ denotes the restricted dual (so for $\mathfrak a$ a finite dimensional Lie algebra, 
$U(\mathfrak a)^\circ$ is the set of linear maps $U(\mathfrak a)\to\mathbb Q$ which vanish on a finite codimensional ideal).  
The assignment $\mathfrak n\mapsto{\sf G}_{\mathfrak n}$ is a functor $\{$pronilpotent $\mathbb Q$-Lie algebras$\}\to\{$affine 
$\mathbb Q$-group schemes$\}$. The Lie algebra of ${\sf G}_{\mathfrak n}$ is then $\mathfrak n$. 
A $\mathbb Q$-group scheme in the image of $\mathfrak n\mapsto{\sf G}_{\mathfrak n}$ 
will be called prounipotent (unipotent if $\mathfrak g$ is nilpotent). 

A projective limit of unipotent $\mathbb Q$-group schemes is prounipotent, and any $\mathbb Q$-group subscheme of a 
prounipotent $\mathbb Q$-group scheme is prounipotent. 

\subsubsection{}\label{sect:lkbacqa}\label{lkbaqa}

Let $\Bbbk$ be a commutative $\mathbb Q$-algebra. One defines $\Bbbk\langle X\rangle$ to be the tensor product algebra 
$\Bbbk\otimes\mathbb Q\langle X\rangle$ and $\Bbbk\langle\langle X\rangle\rangle$ to be the degree completion of 
$\Bbbk\langle X\rangle$, where $x_0$, $x_\sigma$, $\sigma\in\Gamma$ are all of degree 1. The free algebra 
structure on $\Bbbk\langle X\rangle$ induces on $\Bbbk\langle\langle X\rangle\rangle$ the structure of a complete topological algebra. 
The product in this algebra is denoted by $\cdot$. 
The group $(\Bbbk\langle\langle X\rangle\rangle^\times,\cdot)$ of invertible elements of $\Bbbk\langle\langle X\rangle\rangle$ is its subset 
of series whose constant term belongs to the group $\Bbbk^\times$ of invertible elements of $\Bbbk$. 

For $G$ in $\Bbbk\langle\langle X\rangle\rangle^\times$, define ${\rm{aut}}_G$ to be the automorphism of the topological algebra 
$\Bbbk\langle\langle X\rangle\rangle$ given by $x_0\mapsto x_0$, $x_\sigma\mapsto t_\sigma(G)^{-1}x_\sigma t_\sigma(G)$. 

Define a map 
\begin{equation}\label{act}
\Bbbk\langle\langle X\rangle\rangle^\times\times \Bbbk\langle\langle X\rangle\rangle\to
\Bbbk\langle\langle X\rangle\rangle, 
\quad 
(G,H)\mapsto G\circledast H:=G\cdot {\rm{aut}}_G(H).   
\end{equation}
Then one checks that for any $G,H$ in $\Bbbk\langle\langle X\rangle\rangle^\times$ and $K$ in $\Bbbk\langle\langle X\rangle\rangle$, 
one has 
\begin{equation}\label{assoc:circledast}
(G\circledast H)\circledast K=G\circledast(H\circledast K) 
\end{equation}
and that (\ref{act}) restricts and corestricts to a map $(\Bbbk\langle\langle X\rangle\rangle^\times)^2\to
\Bbbk\langle\langle X\rangle\rangle^\times$. 
This implies: 

\begin{lemma} \label{kX:group} (see \cite{Rac}, \S 3.1.2)
$(\Bbbk\langle\langle X\rangle\rangle^\times,\circledast)$ is a group. 
\end{lemma} 

\begin{remark} 
One checks the identity ${\rm{aut}}_{G\circledast H}={\rm{aut}}_G\circ{\rm{aut}}_H$ for any $G,H$ in 
$\Bbbk\langle\langle X\rangle\rangle^\times$. It follows that the map $G\mapsto{\rm{aut}}_G$ defines a morphism 
$(\Bbbk\langle\langle X\rangle\rangle^\times,\circledast)\to{\rm{Aut}}(\Bbbk\langle\langle X\rangle\rangle)^\Gamma$, where 
the target is the group of automorphisms of $\Bbbk\langle\langle X\rangle\rangle$ commuting with the action of $\Gamma$
by $\sigma\mapsto t_\sigma$ (see \S\ref{sect:1:1}). 
\medskip 
\end{remark}

One checks that the functor $\{$commutative $\mathbb Q$-algebras$\}\to\{$groups$\}$, $\Bbbk\mapsto 
(\Bbbk\langle\langle X\rangle\rangle^\times,\circledast)$ is representable by a Hopf algebra, therefore: 

\begin{lemma} \label{kX:group:scheme} 
The functor $\{$commutative $\mathbb Q$-algebras$\}\to\{$groups$\}$, $\Bbbk\mapsto (\Bbbk\langle\langle X\rangle\rangle^\times,\circledast)$ is a 
$\mathbb Q$-group scheme. 
\end{lemma} 

For $\Bbbk$ a commutative $\mathbb Q$-algebra, let $\Bbbk\langle\langle X\rangle\rangle^\times_1$ be the subset of 
$\Bbbk\langle\langle X\rangle\rangle^\times$ of all series with constant term equal to 1. 
Then $\Bbbk\langle\langle X\rangle\rangle^\times_1$ is a subgroup of $(\Bbbk\langle\langle X\rangle\rangle^\times,\cdot)$. 
One checks that $\Bbbk\langle\langle X\rangle\rangle^\times_1$ is a subgroup of $(\Bbbk\langle\langle X\rangle\rangle^\times,
\circledast)$. 

Moreover, the functor $\Bbbk\mapsto\Bbbk\langle\langle X\rangle\rangle^\times_1$ is also represented by a Hopf algebra, therefore

\begin{lemma} The functor $\Bbbk\mapsto(\Bbbk\langle\langle X\rangle\rangle^\times_1,\circledast)$ is a $\mathbb Q$-group 
subscheme of the group scheme $\Bbbk\mapsto(\Bbbk\langle\langle X\rangle\rangle^\times,\circledast)$. 
\end{lemma}

This group scheme is denoted $\sf{MT}$ in \cite{Rac}, so ${\sf{MT}}(\Bbbk)=(\Bbbk\langle\langle X\rangle\rangle^\times_1,\circledast)$. 

For $n\geq 1$, let $\Bbbk\langle\langle X\rangle\rangle^\times_{1,\geq n}$ be the subset of $\Bbbk\langle\langle X\rangle\rangle^\times_1$ 
consisting of all the series which are $\equiv1$ modulo degree $\geq n$. Then $\Bbbk\langle\langle X\rangle\rangle^\times_{1,\geq n}$ is 
a normal subgroup of $(\Bbbk\langle\langle X\rangle\rangle^\times_1,\circledast)$, and the functor $\Bbbk\mapsto(\Bbbk\langle\langle X\rangle\rangle^\times_1/\Bbbk\langle\langle X\rangle\rangle^\times_{1,\geq n},\circledast)$ is a unipotent $\mathbb Q$-group
scheme. It follows: 
\begin{lemma}
The group $\mathbb Q$-scheme ${\sf{MT}}:=(\Bbbk\mapsto(\Bbbk\langle\langle X\rangle\rangle^\times_1,\circledast))$ is prounipotent. 
\end{lemma}

Recall that $\Bbbk\langle\langle X\rangle\rangle$ is equipped with a Hopf algebra structure, with product $\cdot$ and coproduct $\Delta$
such that the elements $x_0,x_\sigma$, $\sigma\in\Gamma$, are primitive. Then $(\Bbbk\langle\langle X\rangle\rangle^\times_1,\cdot)$ 
contains as a subgroup the set of group-like elements for this coproduct, namely $\mathcal G(\Bbbk\langle\langle X\rangle\rangle,\Delta)
=\{G\in\Bbbk\langle\langle X\rangle\rangle^\times|\Delta(G)=G\otimes G\}$. 

The exponential map relative to $\cdot$ sets up a bijection 
$$
(\mathrm{exp})_\Bbbk:\Bbbk\langle\langle X\rangle\rangle_0\stackrel{\sim}{\to}
\Bbbk\langle\langle X\rangle\rangle^\times_1, 
$$
where the index 0 means formal series with vanishing constant term; this bijection restricts to a bijection 
$$
(\mathrm{exp})_\Bbbk:\hat{\mathfrak{Lib}}_\Bbbk(X)\stackrel{\sim}{\to}
\mathcal G(\Bbbk\langle\langle X\rangle\rangle,\Delta), 
$$
where $\hat{\mathfrak{Lib}}_\Bbbk(X)\subset\Bbbk\langle\langle X\rangle\rangle_0$ is the degree completion of 
${\mathfrak{Lib}}_\Bbbk(X):=\mathfrak{Lib}(X)\otimes\Bbbk$. One checks that 
$\mathcal G(\Bbbk\langle\langle X\rangle\rangle,\Delta)=\mathrm{exp}(\hat{\mathfrak{Lib}}_\Bbbk(X))$ 
is a subgroup of $(\Bbbk\langle\langle X\rangle\rangle^\times_1,\circledast)$. Therefore:

\begin{lemma}
The functor $\Bbbk\mapsto(\mathrm{exp}(\hat{\mathfrak{Lib}}_\Bbbk(X)),\circledast)$ is a $\mathbb Q$-group subscheme of 
${\sf{MT}}=(\Bbbk\mapsto(\Bbbk\langle\langle X\rangle\rangle^\times_1,\circledast))$. 
\end{lemma}

As the $\mathbb Q$-group scheme $\sf{MT}$ is prounipotent, we get: 
\begin{lemma}
The $\mathbb Q$-group scheme $\Bbbk\mapsto(\mathrm{exp}(\hat{\mathfrak{Lib}}_\Bbbk(X)),\circledast)$ is prounipotent. 
\end{lemma}

\subsubsection{}\label{sect:4.1.4} The Lie algebras of the $\mathbb Q$-group schemes ${\sf{MT}}=(\Bbbk\mapsto(\Bbbk\langle\langle X\rangle\rangle^\times_1,
\circledast))$ and $\Bbbk\mapsto(\mathrm{exp}(\hat{\mathfrak{Lib}}_\Bbbk(X)),\circledast)$ may be computed using the ring of dual numbers. 
One gets
\begin{equation}\label{Lie:MT}
\mathrm{Lie}({\sf{MT}})=(\mathbb Q\langle\langle X\rangle\rangle_0,\langle,\rangle)=\widehat{\mt}
\end{equation}
and
$$
\mathrm{Lie}(\Bbbk\mapsto(\mathrm{exp}(\hat{\mathfrak{Lib}}_\Bbbk(X)),\circledast))
=(\hat{\mathfrak{Lib}}(X),\langle,\rangle). 
$$
These are the degree completions of the Lie algebras from (\ref{LA:incl}). 

Since the $\mathbb Q$-group schemes ${\sf{MT}}=(\Bbbk\mapsto(\Bbbk\langle\langle X\rangle\rangle^\times_1,
\circledast))$ and $\Bbbk\mapsto(\mathrm{exp}(\hat{\mathfrak{Lib}}_\Bbbk(X)),\circledast)$ are both prounipotent, 
they are isomorphic to the $\mathbb Q$-group schemes ${\sf G}_{\widehat{\mathfrak{mt}}}
=(\Bbbk\mapsto(\Bbbk\langle\langle X\rangle\rangle_0,\mathrm{cbh}_\circledast))$ and 
${\sf G}_{(\hat{\LL}(X),\langle,\rangle)}=(\Bbbk\mapsto(\hat{\mathfrak{Lib}}_\Bbbk(X),\mathrm{cbh}_\circledast))$ 
respectively, where in both cases, $\mathrm{cbh}_\circledast$ 
denotes the CBH series relative to the bracket $\langle,\rangle$. 

\subsubsection{}\label{tQgsib}

The $\mathbb Q$-group scheme isomorphism between ${\sf G}_{\widehat{\mathfrak{mt}}}
$ 
and ${\sf{MT}}
$  is an assignment
$\Bbbk\mapsto (\mathrm{exp}_\circledast)_{\Bbbk}$, where $(\mathrm{exp}_\circledast)_{\Bbbk}$ is a bijection 
$\Bbbk\langle\langle X\rangle\rangle_0\to \Bbbk\langle\langle X\rangle\rangle_1^\times$, 
which intertwines $\mathrm{cbh}_\circledast$ and 
$\circledast$. One checks that the permutation $(\mathrm{exp})_\Bbbk^{-1}\circ(\mathrm{exp}_\circledast)_{\Bbbk}$ of 
$\Bbbk\langle\langle X\rangle\rangle_0$ 
is given by $x\mapsto x+\sum_{n\geq 2}e_n(x)$, where $e_n(x)$ is a homogeneous polynomial self-map of $\Bbbk\langle\langle 
X\rangle\rangle_0$ of degree $n$ with rational coefficients. 

In particular, the differentials at 0 of the two maps $\mathrm{exp}_{\Bbbk},(\mathrm{exp}_\circledast)_{\Bbbk}:\Bbbk\langle\langle X
\rangle\rangle_0\to\Bbbk\langle\langle X\rangle\rangle^\times_1$ coincide.

\subsection{A $\mathbb Q$-group scheme morphism to ${\sf MT}$}\label{sect:aqgsmtm}

\subsubsection{From characters $(\mathrm{exp}(\hat{\mathfrak{Lib}}_{\Bbbk}(X)),\circledast)\to({\Bbbk},+)$ to  
group homomorphisms $(\mathrm{exp}(\hat{\mathfrak{Lib}}_{\Bbbk}(X)),\circledast)$ 
$\to(\Bbbk\langle\langle X\rangle\rangle^\times_1,\circledast)$}

Let $G,H\in\Bbbk\langle\langle X\rangle\rangle^\times$ and $f_G,f_H\in\Bbbk[[x_1]]^\times$, 
$g_G,g_H\in\Bbbk[[x_0]]^\times$. One has 
\begin{align}\label{interm:eq}
& \nonumber (f_G\cdot G\cdot g_G)\circledast(f_H\cdot H\cdot g_H)= 
f_G\cdot G\cdot g_G\cdot \mathrm{aut}_{f_G\cdot G\cdot g_G}(f_H\cdot H\cdot g_H)\\ 
&=f_G\cdot G\cdot g_G\cdot
\mathrm{aut}_{f_G\cdot G\cdot g_G}(f_H)
\mathrm{aut}_{f_G\cdot G\cdot g_G}(H)
\mathrm{aut}_{f_G\cdot G\cdot g_G}(g_H). 
\end{align}
Since 
$$\mathrm{aut}_{f_G\cdot G\cdot g_G}(f_H)=(f_G\cdot G\cdot g_G)^{-1}\cdot f_H\cdot f_G\cdot G\cdot g_G,
$$ 
$$
\mathrm{aut}_{f_G\cdot G\cdot g_G}(H)=\mathrm{aut}_{G\cdot g_G}(H)=(g_G)^{-1}\cdot\mathrm{aut}_G(H)\cdot g_G, \quad
\mathrm{aut}_{f_G\cdot G\cdot g_G}(g_H)=g_H, 
$$
the right hand side of (\ref{interm:eq}) is equal to $f_H\cdot f_G\cdot G\cdot \mathrm{aut}_G(H)\cdot g_G\cdot g_H$. 

Therefore: 
\begin{lemma}\label{lemma:fggfhg}
Let $G,H\in\Bbbk\langle\langle X\rangle\rangle^\times$ and $f_G,f_H\in\Bbbk[[x_1]]^\times$, 
$g_G,g_H\in\Bbbk[[x_0]]^\times$. Then: 
$$
 (f_G\cdot G\cdot g_G)\circledast(f_H\cdot H\cdot g_H)=f_H\cdot f_G\cdot(G\circledast H)\cdot g_G\cdot g_H. 
$$
\end{lemma}

Let $\underline f{}_n,\underline g{}_n:(\mathrm{exp}(\hat{\Lib}_\Bbbk(X)),\circledast)\to(\Bbbk,+)$ ($n\geq 1$) be two collections of group 
homomorphisms. We construct from them group homomorphisms 
$\underline f,\underline g:(\mathrm{exp}(\hat{\Lib}_\Bbbk(X)),\circledast)\to\Bbbk[[t]]^\times$, 
via 
\begin{equation}\label{def:f:g}
\underline f(G):=\mathrm{exp}(\sum_{n\geq 1}\underline f{}_n(G)t^n), 
\quad
\underline g(G):=\mathrm{exp}(\sum_{n\geq 1}\underline g{}_n(G)t^n). 
\end{equation}

Lemma \ref{lemma:fggfhg} implies: 

\begin{lemma}\label{lemma:1:4} If $\underline f{}_n,\underline g{}_n:(\mathrm{exp}(\hat{\Lib}_\Bbbk(X)),\circledast)\to(\Bbbk,+)$ ($n\geq 1$) are
two collections of group homomorphisms and the maps $\underline f,\underline g$ are defined by (\ref{def:f:g}), then the map 
\begin{equation}\label{form:Theta}
G\mapsto \underline f(G)_{|t\to x_1}\cdot G\cdot \underline g(G)_{|t\to x_0}
\end{equation}
defines a group homomorphism 
$$
\Theta_{(\underline{f}{}_n)_n,(\underline g{}_n)_n}:(\mathrm{exp}(\hat{\Lib}_\Bbbk(X)),\circledast)\to(\Bbbk\langle\langle X\rangle\rangle^\times_1,\circledast). 
$$ 
\end{lemma}

Here the indices $|t\to x_i$, $i=0,1$ stand for the homomorphisms $\Bbbk[[t]]\to\Bbbk\langle\langle X\rangle\rangle$, $t\mapsto x_i$ 
($i=0,1$). 

\subsubsection{Construction of a $\mathbb Q$-group scheme morphism $(\Bbbk\mapsto(\mathrm{exp}(\hat{\Lib}_\Bbbk(X)),\circledast))
\to(\Bbbk\mapsto(\Bbbk\langle\langle X\rangle\rangle^\times_1,\circledast))$}\label{fqgsm}

Recall that the additive group is the $\mathbb Q$-group scheme given by $\Bbbk\mapsto(\Bbbk,+)$. Assume that 
$$
{\sf{f}}_n,{\sf{g}}_n:(\Bbbk\mapsto(\mathrm{exp}(\hat{\Lib}_\Bbbk(X)),\circledast))\to\mathbb G_a
$$ 
are $\mathbb Q$-group scheme morphisms; they give rise to compatible families of group homomorphisms 
${\sf{f}}_n(\Bbbk),{\sf{g}}_n(\Bbbk):(\mathrm{exp}(\hat{\Lib}_\Bbbk(X)),\circledast)\to(\Bbbk,+)$ and therefore to 
a group homomorphism 
$$
\Theta_{({\sf{f}}_n(\Bbbk))_n,({\sf{g}}_n(\Bbbk))_n}:(\mathrm{exp}(\hat{\Lib}_\Bbbk(X)),\circledast)\to(\Bbbk\langle\langle X\rangle\rangle^\times_1,\circledast). 
$$
For $g\in\mathrm{exp}(\hat{\Lib}_\Bbbk(X))$, one checks that $\Theta_{({\sf{f}}_n(\Bbbk))_n,({\sf{g}}_n(\Bbbk))_n}(g)$
can be expressed algebraically in terms of the coefficients of $g$. It follows: 
\begin{lemma}\label{tcomtafaqgm}
Let ${\sf{f}}_n,{\sf{g}}_n:(\Bbbk\mapsto(\mathrm{exp}(\hat{\Lib}_\Bbbk(X)),\circledast))\to\mathbb G_a$
be $\mathbb Q$-group scheme morphisms. There exists a $\mathbb Q$-group scheme morphism
$$
\Theta_{({\sf{f}}_n)_n,({\sf{g}}_n)_n}:(\Bbbk\mapsto(\mathrm{exp}(\hat{\Lib}_\Bbbk(X)),\circledast))\to(\Bbbk\mapsto(\Bbbk\langle\langle X\rangle\rangle^\times_1,\circledast)), 
$$
such that for any commutative $\mathbb Q$-algebra $\Bbbk$, the specialization to $\Bbbk$ of $\Theta_{({\sf{f}}_n)_n,({\sf{g}}_n)_n}$ coincides with 
$\Theta_{({\sf{f}}_n(\Bbbk))_n,({\sf{g}}_n(\Bbbk))_n}$. 
\end{lemma}

\begin{remark}
When the ${\sf f}_n,{\sf g}_n$ are trivial, then $\Theta_{({\sf{f}}_n(\Bbbk))_n,({\sf{g}}_n(\Bbbk))_n}$ is nothing but the inclusion 
$(\mathrm{exp}(\hat{\Lib}_\Bbbk(X)),\circledast)\subset(\Bbbk\langle\langle X\rangle\rangle^\times_1,\circledast)$. 
\medskip 
\end{remark}

The Lie algebra homomorphisms associated with ${\sf f}_n,{\sf g}_n$ are characters 
$$
\mathrm{Lie}({\sf f}_n),\mathrm{Lie}({\sf g}_n):(\hat{\Lib}(X),\langle,\rangle)\to\mathbb Q. 
$$
Then (\ref{form:Theta}) implies that the Lie algebra homomorphism associated to $\Theta_{({\sf{f}}_n)_n,({\sf{g}}_n)_n}$ is
$$
\mathrm{Lie}(\Theta_{({\sf{f}}_n)_n,({\sf{g}}_n)_n}):(\hat{\Lib}(X),\langle,\rangle)\to(\mathbb Q\langle\langle X\rangle\rangle_0,
\langle,\rangle), 
$$
\begin{equation}\label{Lie:Theta}
\psi\mapsto \psi+\sum_{n\geq 1}\mathrm{Lie}({\sf f}_n)(\psi)x_1^n+\sum_{n\geq 1}\mathrm{Lie}({\sf g}_n)(\psi)x_0^n.
\end{equation}

\subsection{Construction of a $\mathbb Q$-group scheme morphism $\Theta$ exponentiating $\theta$}\label{sect:coaqgsmtet}

\subsubsection{Exponentiation of some characters $(\hat{\Lib}(X),\langle,\rangle)\to\mathbb Q$}

\begin{lemma}\label{lemma:characters}
Let $n\geq 1$. There are $\mathbb Q$-group scheme morphisms 
$$
{\sf f}_n^0,{\sf g}_n^0:(\Bbbk\mapsto(\mathrm{exp}(\hat{\Lib}_\Bbbk(X)),\circledast))\to\mathbb G_a,
$$ such that for any $\Bbbk$, the group homomorphisms 
${\sf f}_n^0(\Bbbk),{\sf g}_n^0(\Bbbk):(\mathrm{exp}(\hat{\Lib}_\Bbbk(X)),\circledast)\to(\Bbbk,+)$ are given by 
$$
{\sf f}_n^0(\Bbbk):G\mapsto (G|x_0^{n-1}x_1),\quad{\sf g}_1^0(\Bbbk):G\mapsto (G|x_0), \quad {\sf g}_n^0(\Bbbk):G\mapsto
0 \quad \text{for}\quad n\geq 2. 
$$
The associated Lie algebra homomorphisms coincide with the characters $f_n^0,g_n^0:(\hat{\Lib}(X),\langle,\rangle)\to\mathbb Q$ 
defined in Lemma \ref{lemma:tmfglqdb}. 
\end{lemma} 

{\em Proof.} If $\mathfrak G$ is a pronilpotent $\mathbb Q$-Lie algebra, and if $\mathfrak G\to\mathbb Q$
is a Lie algebra homomorphism (character), then there is a unique morphism of $\mathbb Q$-group schemes
$(\Bbbk\to (\mathfrak G\hat\otimes\Bbbk,\mathrm{cbh}))\to\mathbb G_a$ with differential the initial map $\mathfrak G\to\mathbb Q$. 
{\it The homomorphisms $f_n^0,g_n^0$ from Lemma \ref{lemma:tmfglqdb} therefore give rise to 
$\mathbb Q$-group scheme morphisms 
${\sf f}_n^0,{\sf g}_n^0:(\Bbbk\mapsto(\mathrm{exp}(\hat{\Lib}_\Bbbk(X)),\circledast))\to\mathbb G_a$.} 

Let $n\geq 1$ and for $\Bbbk$ a commutative $\mathbb Q$-algebra, let $(\varphi_n)_{\Bbbk}:(\mathrm{exp}(\hat{\Lib}_\Bbbk(X)),\circledast)
\to\Bbbk$ be given by $(\varphi_n)_{\Bbbk}(G):=(G|x_0^{n-1}x_1)$ for $G\in\mathrm{exp}(\hat{\Lib}_\Bbbk(X))$. 

Let $G_0,H\in\mathrm{exp}(\hat{\Lib}_\Bbbk(X))$, where $(\mathrm{log}(G_0)|x_0)=0$. Then 
$$
(G_0\circledast H|x_0^{n-1}x_1)=(G_0|x_0^{n-1}x_1)
+\sum_{k=0}^{n-1}(G_0|x_0^k)(\mathrm{aut}_{G_0}(H)|x_0^{n-1-k}x_1).
$$ 
All the $(G_0|x_0^k)$ are zero for $k>0$, so 
$$
(G_0\circledast H|x_0^{n-1}x_1)=(G_0|x_0^{n-1}x_1)+(\mathrm{aut}_{G_0}(H)|x_0^{n-1}x_1). 
$$
One also computes $(\mathrm{aut}_{G_0}(H)|x_0^{n-1}x_1)=(H|x_0^{n-1}x_1)$, so 
\begin{equation}\label{last:equality}
(G_0\circledast H|x_0^{n-1}x_1)=(G_0|x_0^{n-1}x_1)+(H|x_0^{n-1}x_1). 
\end{equation}
Let now $G\in\mathrm{exp}(\hat{\Lib}_\Bbbk(X))$ be arbitrary. Let $\alpha:=(\mathrm{log}(G_0)|x_0)$. 
Set $G_0:=G\cdot\mathrm{exp}(-\alpha x_0)$. We have $G=G_0\cdot \mathrm{exp}(\alpha x_0)=G_0\circledast 
\mathrm{exp}(\alpha x_0)$. Then 
$$
(G\circledast H|x_0^{n-1}x_1)=(G_0\circledast\mathrm{exp}(\alpha x_0)\circledast H|x_0^{n-1}x_1)
=(G_0|x_0^{n-1}x_1)+(\mathrm{exp}(\alpha x_0)\circledast H|x_0^{n-1}x_1)
$$
after (\ref{last:equality}).  One has $\mathrm{exp}(\alpha x_0)\circledast H=H\cdot \mathrm{exp}(\alpha x_0)$ (recall that 
$x_0$ is central), so this is 
$$
(G_0|x_0^{n-1}x_1)+(H\cdot \mathrm{exp}(\alpha x_0)|x_0^{n-1}x_1). 
$$
The last term is obviously equal to $(H|x_0^{n-1}x_1)$, so 
$$
(G\circledast H|x_0^{n-1}x_1)=(G_0|x_0^{n-1}x_1)+(H|x_0^{n-1}x_1). 
$$
Also, $(G_0|x_0^{n-1}x_1)=(G_0\cdot \mathrm{exp}(\alpha x_0)|x_0^{n-1}x_1)=(G|x_0^{n-1}x_1)$. So finally
$$
(\varphi_n)_{\Bbbk}(G\circledast H)=(G\circledast H|x_0^{n-1}x_1)=(G|x_0^{n-1}x_1)+(H|x_0^{n-1}x_1)=(\varphi_n)_{\Bbbk}(G)+(\varphi_n)_{\Bbbk}(H)
$$
so that {\it $(\varphi_n)_{\Bbbk}$ is a character.}

Let also $(\gamma_n)_{\Bbbk}:(\mathrm{exp}(\hat{\Lib}_\Bbbk(X)),\circledast)
\to\Bbbk$ be given by $(\gamma_1)_{\Bbbk}(G):=(G|x_0)$ for $G\in\mathrm{exp}(\hat{\Lib}_\Bbbk(X))$, and by 
$(\gamma_n)_{\Bbbk}(G):=0$ for $n>1$ and $G\in\mathrm{exp}(\hat{\Lib}_\Bbbk(X))$. 

Let $G,H\in\mathrm{exp}(\hat{\Lib}_\Bbbk(X))$. One has $G\circledast H=G\cdot\mathrm{aut}_G(H)$, and the series $G,\mathrm{aut}_G(H)$ 
have constant term 1, so $(G\circledast H|x_0)=(G|x_0)+(\mathrm{aut}_G(H)|x_0)$. One checks that $(\mathrm{aut}_G(H)|x_0)=(H|x_0)$, so 
$(G\circledast H|x_0)=(G|x_0)+(H|x_0)$. This implies that {\it $(\gamma_1)_{\Bbbk}$ is a character}. It is also clear that {\it the maps 
$(\gamma_n)_{\Bbbk}$ are characters for $n>1$.} 

The assignments $\Bbbk\mapsto(\varphi_n)_\Bbbk$, $\Bbbk\mapsto(\gamma_n)_\Bbbk$ are functorial. The Lie algebra of the 
assignment $\Bbbk\mapsto(\mathrm{exp}(\hat{\Lib}_\Bbbk(X)),\circledast)$ coincides with $(\hat{\Lib}(X),\langle,\rangle)$.

According to \S\ref{tQgsib}, the maps $(\mathrm{exp})_\Bbbk$ and $(\mathrm{exp}_\circledast)_\Bbbk$ have the same differential at origin. 
This enables one to compute {\it the Lie algebra homomorphisms induced by $\Bbbk\mapsto(\varphi_n)_{\Bbbk},(\gamma_n)_\Bbbk$, which 
are found to be equal to $f_n^0,g_n^0$. } 

Since the assignments $\Bbbk\mapsto\varphi_n(\Bbbk)$ and $\Bbbk\mapsto{\sf f}_n^0(\Bbbk)$ on the one hand, 
$\Bbbk\mapsto\gamma_n(\Bbbk)$ and $\Bbbk\mapsto{\sf g}_n^0(\Bbbk)$ on the other hand, are both functorial and induce the same 
Lie algebra homomorphisms, they coincide. 
\hfill\qed\medskip 

\subsubsection{Construction of $\Theta$}\label{sect:cot}

In Lemma \ref{lemma:characters}, $\mathbb Q$-group scheme morphisms ${\sf f}_n^0$, ${\sf g}_n^0$ are constructed. 
Plugging them in the construction of Lemma \ref{tcomtafaqgm}, one obtains a $\mathbb Q$-group scheme morphism
$\Theta_{({\sf f}_n^0)_n,({\sf g}_n^0)_n}:(\Bbbk\mapsto(\mathrm{exp}(\hat{\Lib}_\Bbbk(X)),\circledast))\to(\Bbbk\mapsto(\Bbbk\langle\langle X\rangle\rangle^\times_1,\circledast))$, which will be denoted $\Theta$. According to \S\ref{fqgsm}, the infinitesimal of $\Theta$
is the Lie algebra homomorphism $\theta$ from Prop. \ref{tmtltq}. Summarizing: 

\begin{prop}\label{lemma:constr:Theta}
There exists a $\mathbb Q$-group scheme morphism
$$
\Theta:(\Bbbk\mapsto(\mathrm{exp}(\hat{\Lib}_\Bbbk(X)),\circledast))\to(\Bbbk\mapsto(\Bbbk\langle\langle X\rangle\rangle^\times_1,\circledast)), 
$$
such that for any commutative $\mathbb Q$-algebra $\Bbbk$, the group homomorphism 
$$
\Theta_\Bbbk:(\mathrm{exp}(\hat{\Lib}_\Bbbk(X)),\circledast)\to
(\Bbbk\langle\langle X\rangle\rangle^\times_1,\circledast)
$$ is given by 
\begin{equation}\label{def:Theta:k}
G\mapsto \mathrm{exp}\Big(\sum_{n\geq 1}{(-1)^{n-1}\over n}(G|x_0^{n-1}x_1)x_1^n\Big)
\cdot G\cdot \mathrm{exp}\big(-(G|x_0)x_0\big).
\end{equation}
The differential of $\Theta$ coincides with the Lie algebra homomorphism $\theta$ in Prop. \ref{tmtltq}.
\end{prop}

\section{Comparison of $\mathbb Q$-group schemes}\label{sect:coqgs}

The group analogue of Thm. \ref{main:thm} is Thm. \ref{thm:tsdcwtsoted}, which states that for $\Bbbk$ any commutative 
$\mathbb Q$-algebra, the subgroup $\widetilde{\sf{DMR}}_0(\Bbbk)$ of $(\mathrm{exp}(\hat{\LL}_\Bbbk(X)),\circledast)$ coincides 
with the stabilizer of $\Delta_*$, viewed as an element of the pull-back under $\Theta_\Bbbk$ of the $\sf MT(\Bbbk)$-module  
$\mathrm{Hom}_{\Bbbk}^{\mathrm{cont}}(\Bbbk\langle\langle Y\rangle\rangle,\Bbbk\langle\langle Y\rangle\rangle^{\hat\otimes 2})$. 
This result is proved in several steps. In \S\ref{sect:sgfas}, we show that in the context of certain modules over a pronilpotent Lie algebra, 
stabilizer group functors coincide with the group functors arising from stabilizer Lie algebras. In \S\ref{tgsaimd}, we introduce the modification 
$\widetilde{{\sf DMR}}_0$ of ${\sf DMR}_0$ corresponding to the product of the Lie algebra $\mathfrak{dmr}_0$ with the abelian 
Lie algebra $\mathbb Q x_0\oplus\mathbb Q x_1$. In \S\ref{sect:tmotla}, we construct a module over $(\mathrm{exp}(\hat{\LL}_\Bbbk(X)),
\circledast)$, together with a submodule, to which the results of \S\ref{sect:sgfas} can be applied. In \S\ref{sect:iofktd}, we put together the 
various constructions to prove Thm. \ref{thm:tsdcwtsoted}. 

\subsection{Stabilizer group functors and schemes}\label{sect:sgfas}

\subsubsection{Module functors associated to a module over a pronilpotent $\mathbb Q$-Lie algebra}\label{mfatamo}

Let $\mathfrak n$ be a pronilpotent $\mathbb Q$-Lie algebra. Let $\mathcal C_{\mathfrak n}$ be the the category 
of topological $\mathfrak n$-modules $M$, which are complete and separated, such that the sequence of subspaces 
$(\mathfrak n^k\cdot M)_{k\geq 0}$ is a neighborhood of origin, and such that for each $k\geq 0$, the quotient 
$\mathfrak n^k\cdot M/\mathfrak n^{k+1}\cdot M$ is finite dimensional. The category $\mathcal C_{\mathfrak n}$
is symmetric and monoidal. 

For $\Bbbk$ a commutative $\mathbb Q$-algebra, we define a $\Bbbk$-module $M\hat\otimes\Bbbk:=\lim_{\leftarrow}
(M/\mathfrak n^k\cdot M)\otimes\Bbbk$. The $\Bbbk$-module $M\hat\otimes\Bbbk$ is then a module 
over the group ${\sf G}_{\mathfrak n}(\Bbbk)=(\mathfrak n\hat\otimes\Bbbk,\mathrm{cbh})$, the action of 
$x\in\mathfrak n\hat\otimes\Bbbk$ on $M\hat\otimes\Bbbk$ being given by the element 
$\sum_{n\geq 0}\rho_{\Bbbk}(x)^n/n!\in\mathrm{Aut}_{\Bbbk}(M\hat\otimes\Bbbk)$, where 
$\rho_{\Bbbk}:\mathfrak n\hat\otimes\Bbbk\to\mathrm{End}_{\Bbbk}(M\hat\otimes\Bbbk)$ is the action 
map derived from the action map $\rho:\mathfrak n\to\mathrm{End}_{\mathbb Q}(M)$ by extension of scalars. 
The assignment $\Bbbk\mapsto M\hat\otimes\Bbbk$ is a module functor over the group functor 
$\Bbbk\mapsto{\sf G}_{\mathfrak n}(\Bbbk)$. 

\subsubsection{Stabilizer group functors}\label{sect:sgf}

Let $\mathfrak n$ be a nilpotent Lie algebra, let $M\in\mathrm{Ob}(\mathcal C_{\mathfrak n})$ and let $v\in M$
be a vector. For $\Bbbk$ a commutative $\mathbb Q$-algebra, set 
\begin{equation}\label{def:stab:v}
\mathrm{Stab}(v)(\Bbbk):=\{g\in{\sf G}_{\mathfrak n}(\Bbbk)|g\cdot v=v\text{ (equality in }M\hat\otimes\Bbbk)\}. 
\end{equation}
The assignment $\Bbbk\mapsto\mathrm{Stab}(v)(\Bbbk)$ is then a subgroup functor of the group functor 
$\Bbbk\mapsto{\sf G}_{\mathfrak n}(\Bbbk)$. 

On the other hand, set $\widehat{\mathfrak{stab}}(v):=\{x\in\mathfrak{n}|x\cdot v=0\}$. Then $\widehat{\mathfrak{stab}}(v)$ is a closed Lie 
subalgebra of $\mathfrak{n}$, which is therefore pronilpotent. With this Lie algebra is associated a $\mathbb Q$-group scheme
${\sf G}_{\widehat{\mathfrak{stab}}(v)}$, which is a group subscheme of ${\sf G}_{\mathfrak n}$.   

\begin{lemma}\label{lemma:tsfbk}
The subgroup functors $\Bbbk\mapsto\mathrm{Stab}(v)(\Bbbk)$ and $\Bbbk\mapsto{\sf G}_{\widehat{\mathfrak{stab}}}(\Bbbk)$ 
of the group functor $\Bbbk\mapsto{\sf G}_{\mathfrak n}(\Bbbk)$ are equal, so for any commutative $\mathbb Q$-algebra 
$\Bbbk$, 
one has
$$
\mathrm{Stab}(v)(\Bbbk)={\sf G}_{\widehat{\mathfrak{stab}}(v)}(\Bbbk)
$$
(equality of subgroups of ${\sf G}_{\mathfrak n}(\Bbbk)$).  
\end{lemma}

{\em Proof.} For $k\geq 0$, set $M^{(k)} :=M/{\mathfrak n}^k\cdot M$ and let $v^{(k)}$ be the image of $v\in M$.  Then $M^{(k)}\in\mathrm{Ob}(\mathcal C_{\mathfrak n})$. One associates with the pair $(M^{(k)}, v^{(k)})$ the following data:
\begin{itemize}
\item the Lie subalgebra $\widehat{\mathfrak{stab}}(v)\subset\mathfrak n$, and therefore the $\mathbb Q$-group scheme inclusion 
${\sf G}_{\widehat{\mathfrak{stab}}(v)}\subset{\sf G}_{\mathfrak n}$; 
\item the subgroup functor $\Bbbk\mapsto\mathrm{Stab}(v^{(k)})(\Bbbk)$ of the group functor $\Bbbk\mapsto{\sf G}_{{\mathfrak n}}(\Bbbk)$. 
\end{itemize}
One has 
$\bigcap_{k\geq 0}\widehat{\mathfrak{stab}}(v^{(k)})
=\widehat{\mathfrak{stab}}(v)$, 
which implies that for any commutative $\mathbb Q$-algebra $\Bbbk$, one has
\begin{equation}\label{bksu}
\bigcap_{k\geq 0}{\sf G}_{\widehat{\mathfrak{stab}}(v^{(k)})}(\Bbbk)=
{\sf G}_{\widehat{\mathfrak{stab}}(v)}(\Bbbk)
\end{equation}
(equality of subsets of ${\sf G}_{\mathfrak{n}}(\Bbbk)$).

On the other hand, for any commutative $\mathbb Q$-algebra $\Bbbk$, one has
\begin{equation}\label{bkmsv}
\bigcap_{k\geq 0}\mathrm{Stab}(v^{(k)})(\Bbbk) = 
\mathrm{Stab}(v)(\Bbbk) \end{equation}(equality of subsets of ${\sf G}_{\mathfrak{n}}(\Bbbk)$).
Recall that ${\sf G}_{{\mathfrak n}}(\Bbbk) = ({\mathfrak n}\hat\otimes\Bbbk,\mathrm{cbh})$. Let $x\in{\mathfrak n}\hat\otimes\Bbbk$. 
Then: 
\begin{itemize}
\item 
Assume that $x\in{\sf G}_{\widehat{\mathfrak{stab}}(v^{(k)})}(\Bbbk)$. This means that $x\in\widehat{\mathfrak{stab}}(v^{(k)})\hat\otimes\Bbbk$. Then $x\cdot v^{(k)} = \sum_{n\geq 0}(1/n!)\rho_{\Bbbk}^{(k)}(x)^n\cdot v^{(k)}$ (equality in $M^{(k)}\hat\otimes\Bbbk$), where $\rho_{\Bbbk}^{(k)}:
\mathfrak n\hat\otimes\Bbbk\to\mathrm{End}_\Bbbk(M^{(k)}\hat\otimes\Bbbk)$ is the action map of 
$\mathfrak n\hat\otimes\Bbbk$ on $M^{(k)}\otimes\Bbbk$. 
Since $\rho_\Bbbk^{(k)}(v)=0$, one has $x\cdot v^{(k)}=v^{(k)}$, so $x\in\mathrm{Stab}(v^{(k)})(\Bbbk)$. This shows 
${\sf G}_{\widehat{\mathfrak{stab}}(v^{(k)})}(\Bbbk)\subset\mathrm{Stab}(v^{(k)})(\Bbbk)$. 

\item 
Assume now that  $x\in\mathrm{Stab}(v^{(k)})(\Bbbk)$. Then $x\cdot v^{(k)}=v^{(k)}$ ; moreover, for any $\ell\geq 1$, $x\cdots x\cdot 
v^{(k)}=v^{(k)}$ ($l$ factors $x$). The $\ell$-th power of $x$ in ${\sf G}_{\mathfrak n}(\Bbbk)$ corresponds to $\ell\cdot x\in{\mathfrak 
n}\hat\otimes\Bbbk$. So $(\ell\cdot x)\cdot v^{(k)}=v^{(k)}$, therefore $\sum_{n\geq 0}\ell^n(\rho_\Bbbk^{(k)}(x)^n/n !)(v^{(k)})=v^{(k)}$. 
One has $\rho_\Bbbk^{(k)}(x)^n(v^{(k)})=0$ for $n\geq k$, therefore 
$$
\forall \ell\geq 1,\quad \sum_{n=1}^{k-1} \ell^n(\rho_\Bbbk^{(k)}(x)^n/n !)(v^{(k)})=0. 
$$
As the matrix with elements $(\ell^n)_{\ell\in[1,k-1], n\in[1,k-1]}$ is invertible, one gets $\rho_\Bbbk^{(k)}(x)(v^{(k)})=0$, therefore 
$x\in\widehat{\mathfrak{stab}}(v^{(k)})\hat\otimes\Bbbk\simeq {\sf G}_{\widehat{\mathfrak{stab}}(v^{(k)})}(\Bbbk)$. This shows 
$\mathrm{Stab}(v^{(k)})(\Bbbk)\subset{\sf G}_{\widehat{\mathfrak{stab}}(v^{(k)})}(\Bbbk)$. 
\end{itemize}

All this implies the equality 
${\sf G}_{\widehat{\mathfrak{stab}}(v^{(k)})}(\Bbbk)=\mathrm{Stab}(v^{(k)})(\Bbbk)$, 
which combined with (\ref{bksu}) and (\ref{bkmsv}), yields the conclusion. \hfill\qed\medskip 

\subsection{The group scheme ${\sf DMR}_0$ and its modification $\widetilde{{\sf DMR}}_0$}\label{tgsdaimd}\label{tgsaimd}

The Lie algebra of the affine $\mathbb Q$-group scheme ${\sf MT}$ is $(\mathbb Q\langle\langle X\rangle\rangle_0,\langle,\rangle)$, 
which is pronilpotent. By (\ref{Lie:MT}), we have then 
$$
{\sf G}_{(\mathbb Q\langle\langle X\rangle\rangle_0,\langle,\rangle)}\simeq {\sf MT}=(\Bbbk\mapsto(\Bbbk\langle\langle 
X\rangle\rangle_1^\times,\circledast)). 
$$
Let $\mathfrak G$ be a closed Lie subalgebra of $(\mathbb Q\langle\langle X\rangle\rangle_0,\langle,\rangle)$, let ${\sf G}_{\mathfrak G}$ 
be the corresponding $\mathbb Q$-group subscheme of ${\sf G}_{(\mathbb Q\langle\langle X\rangle\rangle_0,\langle,\rangle)}$ and let 
$\sf G$ be the group subscheme of ${\sf MT}$, image of ${\sf G}_{\mathfrak G}$ under the isomorphism 
${\sf G}_{(\mathbb Q\langle\langle X\rangle\rangle_0,\langle,\rangle)}\simeq {\sf MT}$, so that the diagram 
$$
\xymatrix{
{\sf G}_{(\mathbb Q\langle\langle X\rangle\rangle_0,\langle,\rangle)}\ar^\sim[r]& {\sf MT}\\ 
{\sf G}_{\mathfrak G}\ar@^{{(}->}[u]\ar^\sim[r]& \sf G\ar@^{{(}->}[u]}
$$
commutes. 

Then for any commutative $\mathbb Q$-algebra $\Bbbk$, we have a group inclusion 
$$
({\sf G}(\Bbbk),\circledast)\subset({\sf MT}(\Bbbk),\circledast)=(\Bbbk\langle\langle X\rangle\rangle_1^\times,\circledast). 
$$

Let now $\tilde{\mathfrak G}:={\mathfrak G}+(\mathbb Q x_0\oplus\mathbb Q x_1)$. As $x_0$ and $x_1$ are central in 
$(\mathbb Q\langle\langle X\rangle\rangle_0,\langle,\rangle)$, $\tilde{\mathfrak G}$ is a closed Lie subalgebra of 
$(\mathbb Q\langle\langle X\rangle\rangle_0,\langle,\rangle)$. Let $\tilde{\sf G}$ be the corresponding group subscheme of 
${\sf MT}$. 

\begin{lemma}
For any commutative $\mathbb Q$-algebra $\Bbbk$, we have 
\begin{equation}\label{tgbbbk}
\tilde{\sf G}(\Bbbk)=\{e^{\beta x_1}\cdot g\cdot e^{\alpha x_0}|\alpha,\beta\in\Bbbk,g\in{\sf G}(\Bbbk)\}
\end{equation}
(equality of subsets of ${\sf MT}(\Bbbk)$). In the right-hand side of (\ref{tgbbbk}), the product $\cdot$ is that of 
$\Bbbk\langle\langle X\rangle\rangle^\times$ (see \S\ref{sect:lkbacqa}).
\end{lemma}

\proof 
The algebraic group scheme $\sf{MT}$ is pro-unipotent, therefore if $\Bbbk$ is any commutative $\mathbb Q$-algebra, the exponential map 
$\mathrm{exp}_\circledast$ sets up a bijection $\mathrm{exp}_\circledast:\widehat{\mathfrak{mt}}(\Bbbk)\stackrel{\sim}{\to}{\sf{MT}}(\Bbbk)$ between 
the group of $\Bbbk$-points of this group scheme and the Lie algebra $\widehat{\mathfrak{mt}}(\Bbbk)$ of $\Bbbk$-points of its Lie algebra, which is 
equal to 
$\Bbbk\langle\langle X\rangle\rangle_0$. One has for any $\alpha,\beta\in \Bbbk$, 
\begin{equation}\label{ids:star}
\mathrm{exp}_\circledast(\alpha x_0)=\mathrm{exp}(\alpha x_0), \quad 
\mathrm{exp}_\circledast(\beta x_1)=\mathrm{exp}(\beta x_1)
\end{equation}
and for any $g\in{\sf{MT}}(\Bbbk)$,  
$$
\mathrm{exp}_\circledast(\alpha x_0)\circledast g=g\cdot \mathrm{exp}(\alpha x_0),\quad  
g\circledast\mathrm{exp}_\circledast(\beta x_1)=\mathrm{exp}(\beta x_1)\cdot g
$$ 
(identities in ${\sf{MT}}(\Bbbk)=1+\Bbbk\langle\langle X\rangle\rangle_0$); in the right-hand sides, the products are that of 
$\Bbbk\langle\langle X\rangle\rangle$. If now $\alpha,\beta,\alpha',\beta'\in \Bbbk$ and $g,g'\in{\sf{MT}}(\Bbbk)$, 
\begin{align*}
& (e^{\beta x_1}\cdot g\cdot e^{\alpha x_0})\circledast (e^{\beta' x_1}\cdot g'\cdot e^{\alpha' x_0})\\ & 
= \mathrm{exp}_\circledast(\beta x_1)\circledast g\circledast\mathrm{exp}_\circledast(\alpha x_0)\circledast
\mathrm{exp}_\circledast(\beta' x_1)\circledast g'\circledast\mathrm{exp}_\circledast(\alpha' x_0)
\\ & =\mathrm{exp}_\circledast((\beta+\beta') x_1)\circledast g\circledast g'\circledast\mathrm{exp}_\circledast((\alpha+\alpha') x_0)
=e^{(\beta+\beta') x_1}\cdot (g\circledast g')\cdot e^{(\alpha+\alpha') x_0},  
\end{align*}
where the second equality follows from the centrality of $x_0,x_1$ in $\widehat{\mathfrak{mt}}(\Bbbk)$. 

This implies that the assignment $\Bbbk\mapsto$(right-hand side of (\ref{tgbbbk})) is a subgroup functor of $\Bbbk\mapsto
{\sf MT}(\Bbbk)$. The associated Lie algebra is $\tilde{\mathfrak G}$, so this functor coincides with 
$\Bbbk\mapsto\tilde{\sf G}(\Bbbk)$. 
\hfill \qed\medskip 

Assume that there is an embedding $\iota:\Gamma\to\mathbb C^\times$. With this datum and $\lambda\in \Bbbk$ is associated in 
\cite{Rac}, D\'ef.\ 3.2.1, the set ${\sf{DMR}}_\lambda^\iota(\Bbbk)$. When $\lambda=0$, this set does not depend on $\iota$ 
and will be denoted simply ${\sf{DMR}}_0(\Bbbk)$. 
One has: 

\begin{defin}\label{def:3:6} (see \cite{Rac}, D\'ef. 3.2.1) \label{def:dmr0}
The set ${\sf{DMR}}_0(\Bbbk)$ is the set of elements $G$ of $\Bbbk\langle\langle X\rangle\rangle^\times$, such that:
$$
(G|1)=1, \quad (G|x_0)=(G|x_1)=0, \quad (G|x_0x_1)=0, \quad (G|x_\sigma-x_{\sigma^{-1}})=0 
$$ 
for any $\sigma\in\Gamma$, and 
$$
\Delta(G)=G\otimes G, \quad\Delta_*(G_\star)=G_\star\otimes G_\star
$$
(identities in $\Bbbk\langle\langle X\rangle\rangle^{\hat\otimes 2}$ and $\Bbbk\langle\langle Y\rangle\rangle^{\hat\otimes 2}$), where
\begin{equation}\label{def:G:mult:star}
G_\star=\mathrm{exp}\Big(\sum_{n\geq 2}{{(-1)^{n-1}}\over n}(\pi_Y(G)|y_{n,1})y_{1,1}^n\Big)\cdot\mathbf q(\pi_Y(G))
\end{equation}
(an element of $\Bbbk\langle\langle Y\rangle\rangle^\times$).\footnote{In \cite{Rac}, the notation $\psi\mapsto\psi_\star$ is used to 
denote both the `additive' map $\psi\mapsto \psi_*$, where $\psi_*$ is as in (\ref{def:add:psi:star}) and the `multiplicative' 
map $G\mapsto G_\star$, where $G_\star$ is as in 
(\ref{def:G:mult:star}). We made the choice of giving two different notations to these maps in order to avoid confusions.}
\end{defin}
Then: 

\begin{lemma} (see \cite{Rac}, \S\S 3.2.3 and 3.2.8) This assignment $\Bbbk\mapsto{\sf{DMR}}_0(\Bbbk)$ corresponds to 
a subgroup scheme ${\sf{DMR}}_0$ of $\sf{MT}$, with Lie algebra $\mathfrak{dmr}_0$. 
\end{lemma}

We derive from there: 

\begin{lemma}\label{lemma:dmr0:tilde}\label{lemma:tilde:dmr0}
The subgroup functor of $\Bbbk\mapsto\sf{MT}(\Bbbk)$ corresponding to the degree completion of the Lie algebra $\widetilde{\mathfrak{dmr}}_0:=\mathfrak{dmr}_0
\oplus(\mathbb Q x_0\oplus\mathbb Q x_1)$ is $\Bbbk\mapsto\widetilde{{\sf{DMR}}}_0(\Bbbk)$, where  
$$
\widetilde{{\sf{DMR}}}_0(\Bbbk)=\{e^{\beta x_1}\cdot g\cdot e^{\alpha x_0}\ | \ \alpha,\beta\in \Bbbk, \ g\in{{\sf{DMR}}}_0(\Bbbk)\} 
$$
for $\Bbbk$ any commutative $\mathbb Q$-algebra. 
\end{lemma}

\subsection{Topological modules over the Lie algebra $(\hat{\Lib}(X),\langle,\rangle)$}\label{sect:tmotla}

\subsubsection{Graded modules over $(\Lib(X),\langle,\rangle)$}\label{gmol}

Recall from \S \ref{sect:1:1} the Lie algebras $(\Lib(X),\langle,\rangle)$ and $\mt=(\mathbb Q\langle X\rangle_0,\langle,\rangle)$, from \S \ref{sect:fctalah}
the morphism $\theta:(\Lib(X),\langle,\rangle)\to\mt$; from Lemma \ref{lemma:tmqteqpts} the $\mt$-module structure on $\mathbb Q\langle Y\rangle$
given by the map $\mathbb Q\langle X\rangle_0\ni\psi\mapsto s^Y_\psi\in\mathrm{End}(\mathbb Q\langle Y\rangle)$, from \S \ref{sect:stabilizer}
the induced $\mt$-module structure on $\mathrm{Hom}_{\mathbb Q}(\mathbb Q\langle Y\rangle,\mathbb Q\langle Y\rangle^{\otimes 2})$, 
as well as the pull-back $(\Lib(X),\langle,\rangle)$-module structure on the same space induced by $\theta$. All the objects in this construction 
are $\mathbb Z$-graded: the gradings on $\Lib(X)$ and $\mathbb Q\langle X\rangle_0$ are defined by $\mathrm{deg}(x_\gamma)=1$
for $\gamma\in\Gamma\cup\{0\}$, and the grading on $\mathbb Q\langle Y\rangle$ is defined by $\mathrm{deg}(y_{n,\gamma})=n$
for $n\geq 1$, $\gamma\in\Gamma\cup\{0\}$. The nonzero homogeneous components of $\Lib(X)$ (resp., of $\mt$, of 
$\mathbb Q\langle Y\rangle$) correspond to degrees in $\mathbb Z_{>0}$ (resp., $\mathbb Z_{>0}$, $\mathbb Z_{\geq0}$), 
while $\mathrm{Hom}_{\mathbb Q}(\mathbb Q\langle Y\rangle,\mathbb Q\langle 
Y\rangle^{\otimes 2})$ has nonzero components of all degrees in $\mathbb Z$.  

The element $\Delta_*$ of $\mathrm{Hom}_{\mathbb Q}(\mathbb Q\langle Y\rangle,\mathbb Q\langle Y\rangle^{\otimes 2})$ is of degree 0. 
Let $U(\Lib(X),\langle,\rangle)$ be the universal enveloping algebra of $(\Lib(X),\langle,\rangle)$ and let 
\begin{equation}\label{def:M} 
M :=U(\Lib(X),\langle,\rangle)\cdot\Delta_*. 
\end{equation} 
Then $M$ is a graded sub-$(\Lib(X),\langle,\rangle)$-module of $\mathrm{Hom}_{\mathbb Q}(\mathbb Q\langle Y\rangle,\mathbb Q\langle Y\rangle^{\otimes 2})$; 
it has nonzero components only in degrees $\geq 0$. 

\subsubsection{Topological modules over $(\hat{\Lib}(X),\langle,\rangle)$}\label{fmol}

In \S\ref{lkbaqa}, we introduced the degree completion $(\hat{\Lib}(X),\langle,\rangle)$ of $(\Lib(X),\langle,\rangle)$, and in \S\ref{sect:4.1.4}, 
the degree completion $\widehat{\mt}$ of $\mt$.  Being completions of positively graded Lie algebras, both Lie algebras are pronilpotent. 
The morphism $\theta$ from \S \ref{sect:fctalah} extends to a continuous morphism $$\hat\theta : (\hat{\Lib}(X),\langle,\rangle)\to\widehat{\mt} .$$
Let $\mathbb Q\langle\langle Y\rangle\rangle$ be the degree completion of $\mathbb Q\langle Y\rangle$. Then $\mathbb Q\langle\langle 
Y\rangle\rangle$ is a topological module over $\widehat{\mt}$. Let $\mathbb Q\langle\langle Y\rangle\rangle^{\hat\otimes 2}$ be the degree 
completion of $\mathbb Q\langle\langle Y\rangle\rangle^{\otimes 2}$ and let $\mathrm{Hom}_{\mathbb Q}^{\mathrm{cont}}(\mathbb Q\langle\langle 
Y\rangle\rangle,\mathbb Q\langle\langle Y\rangle\rangle^{\hat\otimes 2})$ be the set of continuous $\mathbb Q$-linear maps $\mathbb 
Q\langle\langle Y\rangle\rangle \to \mathbb Q\langle\langle Y\rangle\rangle^{\hat\otimes 2}$. Using the notation $V[i]$ for the 
homogeneous component of degree $i$ of a $\mathbb Z$-graded $\mathbb Q$-vector space $V$, one identifies this set 
with the subspace of 
$\prod_{i,j\geq 0}\mathrm{Hom}_{\mathbb Q}(\mathbb Q\langle Y\rangle[i],\mathbb Q\langle Y\rangle^{\otimes 2}[j])$ 
of all collections $\underline\varphi=(\varphi_{ij})_{i,j\geq 0}$, satisfying: 
\begin{equation}\label{cond:rest}
\text{there exists a map }F_{\underline\varphi}:\mathbb Z_{\geq 0}\to\mathbb Z_{\geq 0}\text{ with }\limm_{n\to+\infty}
F_{\underline\varphi}(n)=+\infty,\text{ such that } 
\varphi_{ij}=0\text{ if }j<F_{\underline\varphi}(i). 
\end{equation}

The inclusion of $\mathrm{Hom}_{\mathbb Q}(\mathbb Q\langle Y\rangle,\mathbb Q\langle Y\rangle^{\otimes 2})$ in this space may then 
be described as follows: this is the subspace of all collections $(\varphi_{ij})_{i,j\geq 0}$, such that there exists an integer $N\geq 0$, 
such that $\varphi_{ij}=0$ for $|i-j|>N$. 

The space $\mathrm{Hom}_{\mathbb Q}^{\mathrm{cont}}(\mathbb Q\langle\langle Y\rangle\rangle,\mathbb Q\langle\langle 
Y\rangle\rangle^{\hat\otimes 2})$ is then equipped with the structure of a topological $\widehat{\mt}$-module, which can 
be pulled back to a topological $(\hat{\Lib}(X),\langle,\rangle)$-module structure through $\hat\theta$. 

Set $\hat M:=\prod_{i\geq 0}M[i]$, where $M$ is the $(\Lib(X),\langle,\rangle)$-module defined in (\ref{def:M}) and the $M[i]$'s are its homogeneous 
components. Then $\hat M$ is a topological $(\hat{\Lib}(X),\langle,\rangle)$-module. 

The canonical inclusion $M\subset\mathrm{Hom}_{\mathbb Q}(\mathbb Q\langle Y\rangle,\mathbb Q\langle Y\rangle^{\otimes 2})$ 
extends to an inclusion
$$
\hat M\subset \mathrm{Hom}_{\mathbb Q}^{\mathrm{cont}}(\mathbb Q\langle\langle 
Y\rangle\rangle,\mathbb Q\langle\langle Y\rangle\rangle^{\hat\otimes 2})
$$
which is compatible with the topological module structures of both sides over $(\hat{\Lib}(X),\langle,\rangle)$. 

Note that $\hat M$ is an object of the category $\mathcal C_{(\hat{\Lib}(X),\langle,\rangle)}$ associated to the pronilpotent Lie algebra 
$(\hat{\Lib}(X),\langle,\rangle)$ in \S\ref{mfatamo}. 

\subsection{Identification of $\Bbbk\mapsto\widetilde{{\sf DMR}}_0(\Bbbk)$ 
with a stabilizer group functor}\label{sect:iofktd}

Recall from \S1.1 the inclusion of graded Lie algebras $(\Lib(X),\langle,\rangle)\subset\mt$, and from Lemma \ref{lemma:dmr0:tilde} 
and \S1.4 the definitions of the graded Lie subalgebras $\widetilde{\dmr}_0$ and $\mathfrak{stab}(\Delta_*)$ in $(\Lib(X),\langle,\rangle)$.

 In this section, we introduce the following notation $\mathfrak{stab}(v,V):=\{x\in\mathfrak g|x\cdot v=0\}$
for any collection of a graded Lie algebra $\mathfrak g$, a graded $\mathfrak g$-module $V$, and a homogeneous element $v$ in $V$. 
Then $\mathfrak{stab}(\Delta_*)=\mathfrak{stab}(\Delta_*,\mathrm{Hom}_{\mathbb Q}(\mathbb Q\langle 
Y\rangle,\mathbb Q\langle Y\rangle^{\otimes 2}))$, the underlying Lie algebra being $(\Lib(X),\langle,\rangle)$. 

According to Thm. \ref{main:thm}, we have $\mathfrak{stab}(\Delta_*,\mathrm{Hom}_{\mathbb Q}(\mathbb Q\langle 
Y\rangle,\mathbb Q\langle Y\rangle^{\otimes 2}))=\widetilde{\dmr}_0$. 

Recall from \S\ref{gmol} the $(\Lib(X),\langle,\rangle)$-submodule $M\subset\mathrm{Hom}_{\mathbb Q}(\mathbb Q\langle 
Y\rangle,\mathbb Q\langle Y\rangle^{\otimes 2})$. As $\Delta_*\in M$, one has $\mathfrak{stab}(\Delta_*,
\mathrm{Hom}_{\mathbb Q}(\mathbb Q\langle Y\rangle,\mathbb Q\langle Y\rangle^{\otimes 2}))=\mathfrak{stab}(\Delta_*,M)$, therefore 
\begin{equation}\label{eq:stab:dmr0}
\mathfrak{stab}(\Delta_*,M)=\widetilde{\dmr}_0
\end{equation}
In \S\ref{lkbaqa}, we introduced the degree completion $(\hat{\Lib}(X),\langle,\rangle)$ of $(\Lib(X),\langle,\rangle)$, 
and in \S\ref{fmol}, the degree completion $\hat M$ of $M$, which is a topological module over $(\hat{\Lib}(X),\langle,\rangle)$. 
According to \S\ref{sect:sgf}, one attaches to these data and to the element $\Delta_*\in\hat M$ the pronilpotent Lie algebra 
$\widehat{\mathfrak{stab}}(\Delta_*)$, which is a closed Lie subalgebra of $(\hat{\Lib}(X),\langle,\rangle)$. It follows from the fact that 
$M$ is graded and $\Delta_*$ is homogeneous that this Lie algebra coincides with the degree completion of $\mathfrak{stab}(\Delta_*,M)$.

Taking degree completions in (\ref{eq:stab:dmr0}), we obtain 
\begin{equation}
\widehat{\mathfrak{stab}}(\Delta_*)=\widetilde{\dmr}_0^\wedge \quad (\text{equality of closed Lie subalgebras of }(\hat{\Lib}(X),\langle,\rangle)),
\end{equation}
where $\widetilde{\dmr}_0^\wedge$ is the degree completion of $\widetilde{\dmr}_0$. We have therefore an equality of pronilpotent Lie 
algebras 
$\widetilde{\dmr}_0^\wedge=\widehat{\mathfrak{stab}}(\Delta_*) \hookrightarrow (\hat{\Lib}(X),\langle,\rangle)
\hookrightarrow\widehat{\mt}$. 
Applying the functor $\mathfrak n\mapsto {\sf G}_{\mathfrak n}$ of \S\ref{rtaqla} and specializing to a commutative $\mathbb Q$-algebra
$\Bbbk$, one obtains the equality of groups 
\begin{equation}\label{equality:gp:schemes}
{\sf G}_{\widetilde{\dmr}_0^\wedge}(\Bbbk)={\sf G}_{\widehat{\mathfrak{stab}}(\Delta_*)}(\Bbbk)
\hookrightarrow{\sf G}_{(\hat{\Lib}(X),\langle,\rangle)}(\Bbbk)\hookrightarrow{\sf G}_{\widehat{\mt}}(\Bbbk). 
\end{equation}
According to Lemma \ref{lemma:tsfbk} and to (\ref{def:stab:v}), one has 
\begin{equation}\label{equality:G:stab:Delta}
{\sf G}_{\widehat{\mathfrak{stab}}(\Delta_*)}(\Bbbk)=\{g\in{\sf G}_{(\hat{\Lib}(X),\langle,\rangle)}(\Bbbk)\ |\ 
g\cdot\Delta_*=\Delta_*\quad (\text{equality in }\hat M\hat\otimes\Bbbk)\}. 
\end{equation}
According to \S \ref{sect:eotlaht}, there is a group isomorphism  ${\sf G}_{\widehat{\mt}}(\Bbbk)\simeq
{\sf MT}(\Bbbk)=(\Bbbk\langle\langle X\rangle\rangle_1^\times,\circledast)$. 
This group isomorphism takes ${\sf G}_{(\hat{\Lib}(X),\langle,\rangle)}(\Bbbk)$ to 
$(\mathrm{exp}(\hat{\Lib}_\Bbbk(X)),\circledast)$ and 
${\sf G}_{\widetilde{\mathfrak{dmr}}_0^\wedge}(\Bbbk)$ to 
$(\widetilde{{\sf DMR}}_0(\Bbbk),\circledast)$ (see Def. \ref{def:dmr0} and Lemma \ref{lemma:tilde:dmr0}).

On the other hand, (\ref{equality:G:stab:Delta}) implies that the image of 
${\sf G}_{\widehat{\mathfrak{stab}}(\Delta_*)}(\Bbbk)$ is
\begin{equation}\label{descr:image}
\{g\in (\mathrm{exp}(\hat{\Lib}_\Bbbk(X)),\circledast)\ |\ 
g\cdot\Delta_*=\Delta_*\quad(\text{equality in }\hat M\hat\otimes\Bbbk)\}. 
\end{equation} 
The space $\Bbbk\langle\langle Y\rangle\rangle$ is a topological module over 
$(\mathrm{exp}(\hat{\Lib}_\Bbbk(X)),\circledast)$. This leads to a 
$(\mathrm{exp}(\hat{\Lib}_\Bbbk(X)),\circledast)$-module structure over 
$\mathrm{Hom}_\Bbbk^{\mathrm{cont}}(\Bbbk\langle\langle Y\rangle\rangle,
\Bbbk\langle\langle Y\rangle\rangle^{\hat\otimes 2})$.
This space can be identified with the subspace of 
$\prod_{i,j\geq 0}\mathrm{Hom}_\Bbbk^{\mathrm{cont}}(\Bbbk\langle Y\rangle[i],
\Bbbk\langle Y\rangle^{\otimes 2}[j])$ 
defined by conditions (\ref{cond:rest}). Then there is a module inclusion
$$
\hat M\hat\otimes\Bbbk\subset 
\mathrm{Hom}_\Bbbk^{\mathrm{cont}}(\Bbbk\langle\langle Y\rangle\rangle,
\Bbbk\langle\langle Y\rangle\rangle^{\hat\otimes 2}). 
$$
It follows that the module $\hat M\hat\otimes\Bbbk$ may be replaced by 
$\mathrm{Hom}_\Bbbk^{\mathrm{cont}}(\Bbbk\langle\langle Y\rangle\rangle,
\Bbbk\langle\langle Y\rangle\rangle^{\hat\otimes 2})$ in the description of the image of 
${\sf G}_{\widehat{\mathfrak{stab}}(\Delta_*)}(\Bbbk)$ given in (\ref{descr:image}).

The image of the equality in diagram (\ref{equality:gp:schemes}) under the isomorphism 
${\sf G}_{\widehat{\mt}}(\Bbbk)\simeq(\Bbbk\langle\langle X\rangle\rangle_1^\times,\circledast)$ 
is then the equality
\begin{equation}\label{quasi:final}
\widetilde{{\sf DMR}}_0(\Bbbk)
=\{g\in(\mathrm{exp}(\hat{\Lib}_\Bbbk(X)),\circledast)\ |\ g\cdot\Delta_*=\Delta_*\quad
(\text{equality in }\mathrm{Hom}_\Bbbk^{\mathrm{cont}}(\Bbbk\langle\langle Y\rangle\rangle,
\Bbbk\langle\langle Y\rangle\rangle^{\hat\otimes 2}))\}. 
\end{equation} 
 According to (\ref{assoc:circledast}), the map 
\begin{equation}\label{module:kX}
\Bbbk\langle\langle X\rangle\rangle^\times_1\times \Bbbk\langle\langle X\rangle\rangle
\to \Bbbk\langle\langle X\rangle\rangle,\quad (g,H)\mapsto g\circledast H=:S_g(H)
\end{equation}
defines on $\Bbbk\langle\langle X\rangle\rangle$ a structure of 
$(\Bbbk\langle\langle X\rangle\rangle^\times_1\times ,\circledast)$-module. 
The induced $\widehat{\mt}$-module structure over $\Bbbk\langle\langle X\rangle\rangle$
is obtained from the $\mt$-module structure on $\mathbb Q\langle X\rangle$ denoted 
$\psi\mapsto s_\psi$ in \S \ref{sect:1:2:new} by tensor product with $\Bbbk$ and completion. 
This implies that (\ref{module:kX}) is the module functor associated to $\psi\mapsto s_\psi$
according to the procedure of \S \ref{mfatamo}. 
As in \S \ref{sect:1:2:new}, one shows that for $g\in\Bbbk\langle\langle X\rangle\rangle_1^\times$, there exists 
a unique automorphism $S^Y_g$ of $\Bbbk\langle\langle Y\rangle\rangle$, such that the diagram 
$$
\xymatrix{
\Bbbk\langle\langle X\rangle\rangle\ar^{S_g}[r]\ar_{\mathbf q\circ\pi_Y}[d] & 
\Bbbk\langle\langle X\rangle\rangle\ar^{\mathbf q\circ\pi_Y}[d]\\
\Bbbk\langle\langle Y\rangle\rangle\ar_{S_g^Y}[r] &\Bbbk\langle\langle Y\rangle\rangle 
}
$$
commutes.  
The map 
\begin{equation}\label{module:kY}
\Bbbk\langle\langle X\rangle\rangle^\times_1 \times 
\Bbbk\langle\langle Y\rangle\rangle\to \Bbbk\langle\langle Y\rangle\rangle, \quad 
(g,K)\mapsto S^Y_g(K)
\end{equation}
defines on $\Bbbk\langle\langle Y\rangle\rangle$ a structure of 
$(\Bbbk\langle\langle X\rangle\rangle^\times_1,\circledast)$-module, and gives 
rise to the module functor associated to $\psi\mapsto s_\psi^Y$.
The map $S^Y_g$ may be described explicitly as follows
$$
S^Y_g=\mathbf q\circ\tilde S^Y_g\circ\mathbf q^{-1}, 
$$  
where $\tilde S_g^Y\in\mathrm{End}_\Bbbk(\Bbbk\langle\langle Y\rangle\rangle)$ is defined by 
\begin{align*}
& \tilde S_g^Y(K) :=\pi_Y(g\circledast \tilde K)\\
\end{align*}
for $K=K((y_{n,\sigma})_{n\geq 1,\sigma\in\Gamma})\in\Bbbk\langle\langle Y\rangle\rangle$ 
and where $\tilde K :=K((x_0^{n-1}x_\sigma)_{n\geq 1,\sigma\in\Gamma})\in\Bbbk\langle\langle 
X\rangle\rangle$. 
According to Prop. \ref{lemma:constr:Theta}, the map $\Theta_\Bbbk$ given by (\ref{def:Theta:k}) defines a 
group morphism $(\mathrm{exp}(\hat{\Lib}_\Bbbk(X)),\circledast)\to
(\Bbbk\langle\langle X\rangle\rangle_1^\times,\circledast)$. It follows that the map 
$$
(\mathrm{exp}(\hat{\Lib}_\Bbbk(X)),\circledast)\times \Bbbk\langle\langle Y\rangle\rangle
\to\Bbbk\langle\langle Y\rangle\rangle,\quad (h,K)\mapsto S^Y_{\Theta_\Bbbk(h)}(K)
$$
defines on $\Bbbk\langle\langle Y\rangle\rangle$ the structure of a 
$(\mathrm{exp}(\hat{\Lib}_\Bbbk(X)),\circledast)$-module, corresponding to the
 $(\Lib(X),\langle,\rangle)$-module structure on $\mathbb Q\langle Y\rangle$ given by 
$\psi\mapsto s^Y_{\theta(\psi)}$. 

The action of $(\mathrm{exp}(\hat{\Lib}_\Bbbk(X)),\circledast)$ on $\mathrm{Hom}_\Bbbk^{\mathrm{cont}}(\Bbbk\langle\langle Y
\rangle\rangle,
\Bbbk\langle\langle Y\rangle\rangle^{\hat\otimes 2}))$ is then given by $g\cdot u:=
(S^Y_{\Theta_\Bbbk}(g))^{\otimes 2}\circ u\circ(S^Y_{\Theta_\Bbbk}(g))^{-1}$. Therefore the right hand side of 
(\ref{quasi:final}) is equal to 
$$
\{g\in(\mathrm{exp}(\hat{\Lib}_\Bbbk(X)),\circledast)\ |\ (S^Y_{\Theta_\Bbbk}(g))^{\otimes 2}\circ\Delta_*=\Delta_*\circ
S^Y_{\Theta_\Bbbk}(g)\quad (\text{equality in }\mathrm{Hom}_\Bbbk^{\mathrm{cont}}(\Bbbk\langle\langle Y\rangle\rangle,
\Bbbk\langle\langle Y\rangle\rangle^{\hat\otimes 2}))\}. 
$$
Combining this with (\ref{quasi:final}), we get the statement of Thm. \ref{thm:tsdcwtsoted}. 

\thanks{
{\it Acknowledgements}.
This work has been supported by 
JSPS KAKENHI (JP15KK0159)
and Daiko Foundation to HF.
}


\end{document}